\documentclass[12pt]{amsart}
\usepackage{amsmath,amsthm,amsfonts,amssymb,amscd}
\usepackage{array,booktabs,enumitem}
\usepackage{url}
\usepackage[colorlinks=true,linkcolor=blue,citecolor=red]{hyperref}

\numberwithin{equation}{section}
\theoremstyle{plain}
\newtheorem{theorem}{Theorem}[section]
\newtheorem{lemma}[theorem]{Lemma}
\newtheorem{proposition}[theorem]{Proposition}
\newtheorem{corollary}[theorem]{Corollary}
\theoremstyle{definition}
\newtheorem{definition}[theorem]{Definition}
\newtheorem{example}[theorem]{Example}
\theoremstyle{remark}
\newtheorem{remark}[theorem]{Remark}
\newcommand{\func}[1]{\operatorname{#1}}

\begin{document}
\title[Complete Gradient Steady Sasaki--Ricci Solitons]{Complete Gradient
Steady Sasaki--Ricci Solitons}
\author{$^{\ast }$Shu-Cheng Chang}
\address{$^{\ast }$Shanghai Institute of Mathematics and Interdisciplinary
	Sciences, Shanghai 200433, China}
\email{scchang@math.ntu.edu.tw}
\author{$^{\ast\ast }$Jingzhi Tie}
\address{$^{\ast\ast }$Department of Mathematics, University of Georgia,
	Athens, GA 30602-7403, U.S.A.}
\email{jtie@uga.edu}
\author{$^{\ast \dag }$Jiuru Zhou}
\address{$^{\ast \dag }$School of Mathematics, Yangzhou University,
	Yangzhou, 225002, China}
\email{zhoujiuru@yzu.edu.cn}
\thanks{$^\ast$Research supported in part by Startup Foundation for Advanced
	Talents of the Shanghai Institute for Mathematics and Interdisciplinary
	Sciences (No.2302-SRFP-2024-0049). $^{\ast \dag }$Partially supported by
	NSFC Grants 12171417.}

\begin{abstract}
We study complete gradient steady Sasaki-Ricci solitons with nonnegative
transverse holomorphic bisectional curvature. We prove that positive
transverse Ricci curvature implies the existence of a transverse critical
point in every dimension. In dimension five, we give a direct construction
of a complete smooth transverse development for every soliton under the
bisectional curvature assumption. The developed metric is flat, a product of a cigar and a flat
complex line, a product of two scaled cigars, or an irreducible steady K\"{a}%
hler surface with positive Ricci curvature. We also classify solitons with
nonnegative transverse sectional curvature and linear scalar curvature decay
measured from a lifted Reeb leaf. The transversally flat solitons are
Heisenberg quotients. In the positive-Ricci developed case, the total Sasakian
soliton is a line contactization or a suspension determined by a
potential-preserving holomorphic isometry. Explicit lifts illustrate
the classification and yield formulas for Reeb-orbit periods and volume
growth.
\end{abstract}

\maketitle
\tableofcontents

\section{Introduction}

\label{sec:introduction}

Sasakian geometry is the odd-dimensional counterpart of K\"{a}hler geometry.
The Reeb foliation of a Sasakian manifold carries a transverse K\"{a}hler
structure, and many questions about Sasakian metrics can be studied through
this structure (see \cite{BoyerGalicki}). Sasaki-Ricci solitons are
self-similar solutions of the Sasaki--Ricci flow and provide natural models
for its limiting geometry. In this paper we study the complete gradient
steady case, with particular emphasis on the relation between transverse
curvature and the global geometry of the Reeb foliation.

The steady K\"{a}hler-Ricci soliton equation has a rich family of complete
solutions. Hamilton's cigar is the basic nonflat surface model. Cao
constructed rotationally symmetric examples on $\mathbb{C}^{n}$ \cite%
{Cao1996} and established a dimension-reduction theorem for the K\"{a}%
hler-Ricci flow \cite{Cao2004}. Under nonnegative sectional curvature,
Deng-Zhu obtained a classification of steady Ricci solitons with linear
scalar curvature decay \cite{DengZhu2020}. More recently, the examples of
Apostolov-Cifarelli \cite{ApostolovCifarelli} showed that the complete
irreducible steady K\"{a}hler metrics form a larger family than the
rotationally symmetric one. These results suggest a structural
classification in which the transverse metric and the global Sasakian
lifting data are treated separately.

The main difficulty is global. Local soliton identities hold in transverse K%
\"{a}hler charts for every Reeb type. A regular Reeb action gives a smooth
complete quotient, while a quasi-regular action gives an orbifold quotient.
An irregular foliation need not have a Hausdorff leaf space. Thus a local
splitting does not by itself provide the complete K\"{a}hler factors needed
for a global classification. Sections~\ref{sec:quasi-regular}
and~\ref{sec:irregular} give a direct construction in dimension five of a
complete smooth transverse development: a universal covering
of the total space admits a submersion onto a complete simply connected K\"{a%
}hler manifold, with connected fibers equal to the lifted Reeb leaves and
with both the transverse metric and potential pulled back from the base. The
precise definition is given in Definition~\ref{def:development}.

The point of view of this paper is to study both the classification of the
complete transverse K\"ahler steady soliton and the global reconstruction of
the Sasakian manifold from its Reeb foliation. Our emphasis is the
dimension-five structural description and the explicit reconstruction data:
monodromy up to conjugacy, Reeb-orbit periods, quasi-regular approximation,
and volume asymptotics. The smooth K\"ahler splitting and scalar-decay
results used below follow from \cite{Cao2004,DengZhu2020}, and biholomorphic
uniformization follows from \cite{ChauTam2005} when the developed
positive-Ricci soliton has a critical point.

The present paper and the companion work \cite{Rigidity2026} share the
arguments for Ricci monotonicity, critical-point existence, and development
in the positive-Ricci case. Lemma~\ref{lem:transverse-ricci-monotonicity} and
Theorem~\ref{thm:positive-ricci-critical} correspond to Lemma~3.1 and
Theorem~3.1 of that work, and the critical-leaf construction in
Proposition~\ref{prop:critical-leaf-development} and the reconstruction in
Theorem~\ref{thm:steady-suspension-classification} overlap with its
Theorem~3.2 under the common hypotheses. We include the proofs to give a
direct dimension-five construction and a complete account of the lifting
data. Proposition~\ref{prop:critical-leaf-development} also retains a
useful distinction: for compact Reeb closure, it assumes positive transverse
Ricci curvature and an existing critical point, but does not require
nonnegative transverse bisectional curvature. The explicit conjugacy
criterion, the monodromy approximation, the orbifold constructions, and the
leading volume coefficients are developed here in
Sections~\ref{sec:quasi-regular}--\ref{sec:asymptotic}.

For the arbitrary-dimensional decay classification we use the following
global-development result from \cite[Theorem~7.2]{Rigidity2026}. If a
gradient steady Sasaki--Ricci soliton is complete and $B^T\geq0$, then its
lifted Reeb flow on the Riemannian universal cover is a free proper
$\mathbb R$-action. The quotient
\[
D:\widetilde M\longrightarrow N=\widetilde M/\mathbb R
\]
is a complete simply connected smooth K\"ahler steady soliton, with the
transverse tensors and the lifted potential pulled back from $N$. Moreover,
the principal bundle is trivial and, in a global fiber coordinate $t$,
\[
\widetilde M\cong N\times\mathbb R,\qquad
\widetilde\eta=dt+D^*\beta,\qquad d\beta=2\omega_N.
\]
This theorem requires neither volume noncollapsing nor scalar curvature
decay. It supplies the development used in Theorem~\ref{thm:intro-decay}
and Corollary~\ref{cor:intro-decay-total}. The local conservation identities
in Theorem~\ref{thm:conservation} are proved in transverse charts and do not
use this global result. The proof of Theorem~\ref{thm:intro-five} in
Sections~\ref{sec:quasi-regular} and~\ref{sec:irregular} gives an independent,
explicit construction in dimension five; the companion's
Proposition~7.1 and Theorem~7.2 establish existence in every dimension by
globally integrating the flat transverse Ricci kernel.

We explain more precisely the relation with the work of Deng--Zhu. Their
classification theorem \cite[Theorem~1.1]{DengZhu2020} concerns a complete
smooth steady Ricci soliton with nonnegative sectional curvature and a
linear scalar curvature upper bound. We apply its K\"ahler form to the
universal transverse development supplied by the preceding result.
Distance to a lifted Reeb leaf becomes distance from a point on this
development, and the resulting decay excludes Euclidean factors in the
nonflat case. We then recover the total space and its Reeb holonomy from
the developed soliton. The smooth K\"ahler critical-point result is due to
Deng--Zhu \cite[Theorem~1.2]{DengZhu2015}; the transverse argument reproduced
here uses Ricci monotonicity and finite-length gradient trajectories before
a global development is available. The classification up to isometry of
the irreducible positive-Ricci K\"ahler factors remains open.

Our first main result combines the transverse classification with the
development construction.

\begin{theorem}
\label{thm:intro-five} Let $(M^{5},g,\eta ,\xi ,\Phi ,f)$ be a complete
gradient steady Sasaki-Ricci soliton with the nonnegative transverse
holomorphic bisectional curvature 
\begin{equation*}
B^{T}\geq 0.
\end{equation*}%
Then it admits a complete smooth transverse development $(N^{4},h,J,f_{N})$.
The developed K\"{a}hler metric is holomorphically isometric to one of the
following :

\begin{enumerate}
\item 
\begin{equation*}
(\mathbb{C}^{2},g_{\mathrm{flat}}),
\end{equation*}

\item 
\begin{equation*}
(\mathbb{C},ag_{\mathrm{cig}})\times (\mathbb{C},g_{\mathrm{flat}}),
\end{equation*}%
for some $a>0$. Here the cigar metric $g_{\mathrm{cig}}=|\mathrm{d}%
z|^{2}/(1+|z|^{2})$,

\item 
\begin{equation*}
(\mathbb{C},ag_{\mathrm{cig}})\times (\mathbb{C},bg_{\mathrm{cig}}),
\end{equation*}%
for some $a,b>0$,

\item a complete simply connected irreducible K\"{a}hler steady surface with
nonnegative holomorphic bisectional curvature and positive Ricci curvature.
\end{enumerate}

If the Reeb foliation is quasi-regular, its effective orbifold quotient is
developable. That is, there are no non-developable quasi-regular gradient
steady Sasaki-Ricci solitons of $B^{T}\geq 0$.
\end{theorem}

Theorem~\ref{thm:intro-five} is a structural classification. The metric in
the irreducible positive-Ricci case remains to be classified.

The existence of a transverse critical point is the main step to the
quasi-regular or irregular foliation in the case of positive transverse
Ricci curvature. This part of the argument holds in every dimension.

\begin{theorem}
\label{thm:intro-positive} Let $(M^{2n+1},g,\eta ,\xi ,\Phi ,f)$ be a
complete gradient steady Sasaki-Ricci soliton with%
\begin{equation*}
B^{T}\geq 0\text{ \ \ \textrm{and} \ }\func{Ric}^{T}>0.
\end{equation*}
Then

\begin{enumerate}
\item $f$ has a transverse critical point, and 
\begin{equation*}
R^{T}+|\nabla ^{T}f|^{2}=C_{0},\qquad C_{0}=\max_{M}R^{T}.
\end{equation*}

\item The soliton has a complete smooth transverse development $%
(N^{2n},h,\omega _{h},J,f_{N})$, where $N$ is biholomorphic to $\mathbb{C}%
^{n}$ and $f_{N}$ has a unique critical point.

\item There is a compact group $K$ of holomorphic isometries of $N$
preserving $f_{N}$, and a $K$-invariant one-form $\beta $ with $d\beta
=2\omega _{h}$, such that $M$ is isomorphic either to the line
contactization $N\times \mathbb{R}$ or to 
\begin{equation*}
M_{A,\ell }=(N\times \mathbb{R})/\langle (x,t)\mapsto (Ax,t+\ell )\rangle
,\qquad A\in K,\quad \ell >0,
\end{equation*}%
with $\eta =dt+\beta $ and $g=h+\eta ^{2}$.
\end{enumerate}
\end{theorem}

We next consider scalar curvature decay measured on the universal cover.
Let $q:(\widetilde M,\widetilde g)\to (M,g)$ be the Riemannian universal
cover and let $\widetilde{\mathcal L}_{\widetilde o}$ be the orbit through
$\widetilde o$ of the lifted Reeb flow. Write
\begin{equation*}
\rho(\widetilde x)=d_{\widetilde g}(\widetilde x,
\widetilde{\mathcal L}_{\widetilde o}),\qquad \widetilde R^{T}=q^*R^T.
\end{equation*}

\begin{theorem}
\label{thm:intro-decay} Let $(M^{2n+1},g,\eta,\xi,\Phi,f)$, $n\geq 1$, be a
connected complete gradient steady Sasaki--Ricci soliton with
$\func{Sec}^{T}\geq 0$. Suppose there are a point $\widetilde o\in\widetilde M$
and constants $C,r_0>0$ such that
\begin{equation}
\widetilde R^T(\widetilde x)\leq \frac{C}{\rho(\widetilde x)}\qquad
(\rho(\widetilde x)\geq r_0).  \label{eq:lifted-leaf-decay}
\end{equation}
Then the lifted Reeb flow is a free proper $\mathbb{R}$-action and
$D:\widetilde M\to N=\widetilde M/\mathbb{R}$ is a complete smooth transverse
development. Its induced gradient steady K\"{a}hler--Ricci soliton
$(N^{2n},h,J,f_N)$ has one of the following forms, up to holomorphic isometry.

\begin{enumerate}
\item $(\mathbb{C}^{n},g_{\mathrm{flat}})$, with $f_N$ affine.

\item If $n=1$, a scaled Hamilton cigar $(\mathbb{C},ag_{\mathrm{cig}})$,
with $a>0$.

\item If $n\geq 2$, a complete simply connected irreducible K\"{a}hler
steady soliton with $\func{Sec}_{h}\geq 0$ and $\func{Ric}_{h}>0$. For
$o=D(\widetilde o)$ there are constants $c,C',r_1>0$ such that
\begin{equation*}
\frac{c}{d_h(o,y)}\leq R_h(y)\leq \frac{C'}{d_h(o,y)}\qquad
(d_h(o,y)\geq r_1).
\end{equation*}
\end{enumerate}

In the nonflat cases $f_N$ has a unique critical point, its global maximum.
The principal $\mathbb{R}$-bundle $D$ is trivial. In a global trivialization
$\widetilde M=N\times\mathbb{R}$, its lifted tensors satisfy
\begin{equation*}
\widetilde\eta=dt+D^*\beta,\qquad d\beta=2\omega_h,\qquad
\widetilde g=D^*h+\widetilde\eta^{\,2},\qquad \widetilde f=D^*f_N.
\end{equation*}
\end{theorem}

\begin{corollary}
\label{cor:intro-decay-total} Under the hypotheses of Theorem~\ref{thm:intro-decay},
the total Sasakian soliton $M$ has the following forms.

\begin{enumerate}[beginpenalty=10000]
\item If $N$ is flat, $M$ is a quotient of the standard Heisenberg
Sasakian manifold $H^{2n+1}=\mathbb{C}^n\times\mathbb{R}$. In coordinates
$(z,t)$, put
\begin{equation*}
\eta_H=dt+\operatorname{Im}\langle z,dz\rangle,\qquad
 g_H=|dz|^2+\eta_H^2,\qquad \xi_H=\partial_t,\qquad
 f_H=\operatorname{Re}\langle a,z\rangle+b,
\end{equation*}
where $\langle w,z\rangle=\sum_j\overline{w_j}z_j$, $a\in\mathbb{C}^n$ and
$b\in\mathbb{R}$. The deck group is a discrete group acting freely and
properly discontinuously by transformations
\begin{equation}
(z,t)\longmapsto
\bigl(Az+v,\,t+s-\operatorname{Im}\langle v,Az\rangle\bigr),
\qquad A\in U(n),\quad v\in\mathbb{C}^n,\quad s\in\mathbb{R},
\label{eq:flat-decay-deck}
\end{equation}
with $Aa=a$ and $\operatorname{Re}\langle a,v\rangle=0$. For constant
potential this includes Heisenberg lattice quotients.

\item If $N$ is nonflat, it is biholomorphic to $\mathbb{C}^n$. Its group
$K$ of potential-preserving holomorphic isometries is compact, and one can
choose $\beta$ with $d\beta=2\omega_h$ invariant under $K$. Then $M$ is
isomorphic to the line contactization $N\times\mathbb{R}$, with proper Reeb
$\mathbb{R}$-action and trivial fundamental group, or to the suspension
\begin{equation*}
M_{A,\ell}=(N\times\mathbb{R})/\langle (y,t)\mapsto (Ay,t+\ell)\rangle,
\qquad A\in K,\quad \ell>0,
\end{equation*}
with fundamental group $\mathbb{Z}$. The tensors are induced by
$\eta=dt+\beta$, $\xi=\partial_t$, $g=h+\eta^2$ and $f=f_N$, with
$\Phi=-J^{\mathrm{hor}}$ on $\ker\eta$ and $\Phi\xi=0$.

In a suspension, the leaf through $[(y,t)]$ is periodic precisely when
$A^jy=y$ for some $j\geq 1$. Its minimal period is $j\ell$ for the least
such $j$. In particular the critical Reeb leaf has period $\ell$. If $A$
is the identity, the flow is regular. If $A$ has finite order $m>1$, it is
quasi-regular with generic period $m\ell$ and effective orbifold quotient
$[N/\langle A\rangle]$. If $A$ has infinite order, the flow is irregular
and its closure is the torus $(C_A\times\mathbb{R})/\langle(A,\ell)\rangle$,
of dimension $1+\dim C_A^0$, where $C_A=\overline{\langle A\rangle}$.

For a fixed transverse soliton, two suspensions are isomorphic as Sasakian
solitons precisely when their positive translation lengths agree and their
monodromies are conjugate in $K$.
\end{enumerate}

Conversely, each of these models, formed with the indicated group data and a
base from Theorem~\ref{thm:intro-decay}, satisfies the hypotheses of that theorem.
\end{corollary}

In dimension five, Theorem~\ref{thm:intro-five} describes the transverse
development under $B^T\geq 0$. The sectional curvature hypothesis of
Theorem~\ref{thm:intro-decay} implies $B^T\geq 0$, and its lifted-leaf
decay excludes both cigar products by Proposition~\ref{prop:product-obstruction}.
For the remaining flat and irreducible developments,
Corollary~\ref{cor:intro-decay-total} gives the global models for $M$ and
determines the Reeb type in the nonflat case.

The proof of Theorem~\ref{thm:intro-decay} is given in
Section~\ref{subsec:quotient-decay}, where we compare decay on the development
with decay on a smooth quotient. Corollary~\ref{cor:intro-decay-total} is
proved in Section~\ref{subsec:suspension-classification}. Even positive
sectional curvature and two-sided linear scalar curvature decay do not
characterize Cao's metric, as Proposition~\ref{prop:ac} shows.

The relation between global development and the remaining
metric-classification problem is as follows.

\begin{enumerate}
\item Complete smooth transverse developments exist in every dimension
under $B^T\geq0$, by
\cite[Proposition~7.1 and Theorem~1.1(2)]{Rigidity2026}.
In the nonflat intermediate-rank case, the flat Ricci kernel is integrated
globally on the simply connected total space, and the complementary
Sasakian factor is complete and simply connected with positive transverse
Ricci curvature. This resolves the higher-dimensional existence question.
The direct dimension-five construction in
Sections~\ref{sec:quasi-regular} and~\ref{sec:irregular} gives an explicit
description of the local orbifold and cigar-line geometry.

\item How can one classify irreducible steady K\"{a}hler factors with
positive Ricci curvature? Proposition~\ref{prop:ac} shows that positive
sectional curvature and two-sided linear scalar decay allow non-Cao metrics.
\end{enumerate}

Section~\ref{sec:geometry} fixes the conventions. Section~\ref{sec:regular}
treats complete smooth quotients and curvature decay on the development.
Sections~\ref{sec:quasi-regular} and~\ref{sec:irregular} construct transverse
developments and determine the
global lifting data. Section~\ref{sec:examples} gives explicit metrics and
lifting constructions, and Section~\ref{sec:asymptotic} studies volume
growth and asymptotic geometry.

\section{Preliminaries}

\label{sec:geometry}

We fix our conventions and collect the transverse identities used below. For
Sasakian background, we refer to \cite[Chapters~7--8]{BoyerGalicki}.

\subsection{Sasakian Geometry}

\label{subsec:sasaki-preliminaries}

Let $(M^{2n+1},g,\eta,\xi,\Phi)$ be a connected Sasakian manifold. We use $%
\Phi=-\nabla^g\xi$, so that 
\begin{gather}
\eta=g(\xi,\cdot),\qquad \Phi^2=-\func{Id}+\eta\otimes\xi,\qquad
d\eta(X,Y)=2g(X,\Phi Y),  \label{eq:sasaki-algebra} \\
(\nabla^g_X\Phi)Y=g(X,Y)\xi-\eta(Y)X.  \label{eq:sasaki-differential}
\end{gather}
The Reeb field is unit, Killing and geodesic. Our transverse conventions are 
\begin{equation}
\mathcal{D}=\ker\eta,\qquad g^T=g-\eta\otimes\eta,\qquad J^T=-\Phi|_{%
\mathcal{D}},\qquad \omega^T=\tfrac12d\eta.
\label{eq:transverse-conventions}
\end{equation}
Here $\omega^T(X,Y)=g^T(J^TX,Y)=-g^T(\Phi X,Y)=g^T(X,\Phi Y)$ for horizontal
vectors, and $g^T$ also denotes the tensor on $TM$ with kernel $\mathbb{R}%
\xi $. We identify vectors and one-forms using the relevant metric. Our
curvature conventions are 
\begin{gather}
\mathcal{R}(X,Y)Z=\nabla_X\nabla_YZ-\nabla_Y\nabla_XZ -\nabla_{[X,Y]}Z, 
\notag \\
\func{Rm}(X,Y,Z,W)=g(\mathcal{R}(X,Y)Z,W),\qquad K(X,Y)=g(\mathcal{R}%
(X,Y)Y,X)  \label{eq:curvature-conventions}
\end{gather}
for orthonormal real $X,Y$. Nonnegative holomorphic bisectional curvature
means $g(\mathcal{R}(X,JX)JY,Y)\geq0$, equivalently $\func{Rm}(Z,\overline
Z,W,\overline W)\geq0$ for $(1,0)$ vectors $Z,W$.

The Reeb foliation $\mathcal{F}_\xi$ is locally defined by submersions onto
K\"ahler manifolds, with local holomorphic isometries as transition maps.
Transverse tensors and differential operators are defined in these charts.
Identifying the normal bundle with $\mathcal{D}$, the transverse connection
is given by 
\begin{equation*}
\nabla^T_XY=\pi_{\mathcal{D}}(\nabla^g_XY),\qquad \nabla^T_\xi Y=\pi_{%
\mathcal{D}}[\xi,Y],
\end{equation*}
where $X,Y$ are horizontal. We write $B^T\geq0$ for nonnegative transverse
bisectional curvature; it implies $\func{Ric}^T\geq0$ and $R^T\geq0$, with $%
R^T=\func{tr}_{g^T}\func{Ric}^T$ the real scalar trace.

A smooth function $u$ is basic when $\xi u=0$. Since the Reeb leaves are
geodesic, the standard transverse formulas give 
\begin{equation}
\nabla^g u=\nabla^T u,\qquad \func{Hess}_g u|_{\mathcal{D}\times\mathcal{D}}=%
\func{Hess}^T u,\qquad \Delta_g u=\Delta^T u,  \label{eq:basic-ambient}
\end{equation}
where $\Delta^T=\func{tr}_{g^T}\func{Hess}^T$ and $\Delta_g=\func{div}%
_g\nabla^g$. The Sasakian Ricci identities \cite[Chapter~7]{BoyerGalicki}
are 
\begin{equation}
\begin{split}
\func{Ric}_g(\xi,\xi)&=2n,\qquad \func{Ric}_g(X,\xi)=0, \\
\func{Ric}_g(X,Y)&=\func{Ric}^T(X,Y)-2g(X,Y),\qquad R_g=R^T-2n
\end{split}
\label{eq:oneill}
\end{equation}
for horizontal $X,Y$. Thus $B^T\geq0$ implies $\func{Ric}_g\geq-2g$. For
basic $f,u$, the transverse drift operator $L_f=\Delta^T-\langle \nabla^T
f,\nabla^T\,\cdot\,\rangle$ agrees on $u$ with $\Delta_g-\langle \nabla^g
f,\nabla^g\,\cdot\,\rangle$, allowing ambient maximum principles even for
irregular foliations.

\subsection{Reeb Foliations and Complete Transverse Developments}

\label{subsec:reeb-preliminaries}

A Sasakian structure is \emph{quasi-regular} if all leaves of $\mathcal{F}%
_{\xi }$ are compact so that $\xi $ generates a locally free circle action.
It is \emph{regular} if this action is free, that is, the Reeb flow is a
free effective circle action with a common minimal period $L>0$. We call the
structure \emph{irregular} if some Reeb leaf is not closed. On a noncompact
manifold there is also the possibility of a free proper $\mathbb{R}$-action.
This case has a smooth global quotient and is treated separately.

Suppose that $g$ is complete, and let $H=\overline{\{\varphi_t:t\in\mathbb{R}%
\}}\subset\func{Isom}(M,g)$ be the closure of the Reeb flow. By continuity,
every continuous basic function is $H$-invariant. The same observation
applies to smooth Reeb-invariant tensors.

The following is well-known, please refer to \cite[Proposition 2.2]{ChangLinQiu2026-2}.
\begin{proposition}
\label{prop:steady-second-structure} Exactly one of the following occurs:

\begin{enumerate}
\item $H\cong \mathbb{S}^{1}$: the Reeb flow is periodic and quasi-regular,

\item $H\cong \mathbb{R}$ : the Reeb action is free and proper, with
complete smooth K\"{a}hler quotient,

\item $H\cong \mathbb{T}^{k}$, $2\leq k\leq n+1$: the flow is irregular and
its orbits are dense in the compact $H$-orbits.
\end{enumerate}

The group $H$ preserves $g,\eta ,\xi ,\Phi ,f$. In the torus cases, put $%
\mathfrak{h}=\func{Lie}(H)$ and $\mu _{\zeta }=\eta (\zeta )$. These
functions are basic and $H$-invariant, with 
\begin{equation}
\iota _{\zeta }d\eta =-d\mu _{\zeta },\qquad \nabla ^{T}\mu _{\zeta }=2\Phi
\zeta _{\mathcal{D}},\qquad |\nabla ^{T}\mu _{\zeta }|=2|\zeta _{\mathcal{D}%
}|.  \label{eq:steady-moment-map}
\end{equation}%
For $\zeta \in \mathfrak{h}^{+}=\{\zeta \in \mathfrak{h}:\inf_{M}\mu _{\zeta
}>0\}$, the type-I Sasakian deformation is 
\begin{equation}
\eta _{\zeta }=\frac{\eta }{\mu _{\zeta }},\quad \xi _{\zeta }=\zeta ,\quad
\Phi _{\zeta }=\Phi -\Phi (\zeta )\otimes \eta _{\zeta },\quad g_{\zeta }=-%
\tfrac{1}{2}d\eta _{\zeta }(\,\cdot \,,\Phi _{\zeta }\,\cdot \,)+\eta
_{\zeta }^{2}.  \label{eq:steady-type-I}
\end{equation}%
If the moment functions of a basis of $\mathfrak{h}$ are bounded, rational
vectors $\zeta _{i}\rightarrow \xi $ can be chosen in $\mathfrak{h}^{+}$,
giving quasi-regular structures converging in $C^{\infty }$ on compact
subsets. If the transverse gradients of the basis moment functions are also
bounded, the metrics $g_{\zeta _{i}}$ are complete. The resulting structures
are Sasakian, but they need not satisfy the steady soliton equation.
\end{proposition}

For a proper line action or a regular circle action, the projection $%
\pi:M\to Q=M/\mathcal{F}_\xi$ is a Riemannian submersion onto a complete
smooth K\"ahler manifold. The basic potential descends and satisfies 
\begin{equation}
\func{Ric}_{h_Q}+\func{Hess}_{h_Q}f_Q=0.  \label{eq:quotient-soliton}
\end{equation}
A quasi-regular quotient is a complete K\"ahler orbifold, whose universal
orbifold cover is smooth only when it is developable. For an irregular flow, 
$M/H$ has generic dimension $2n+1-\dim H$, rather than the transverse
dimension $2n$, and does not replace a K\"ahler leaf space.

\begin{definition}
\label{def:development} A complete transverse development of a Sasakian
soliton consists of a connected covering $q:\widehat{M}\rightarrow M$ and a
surjective submersion $D:\widehat{M}\rightarrow N$ onto a complete simply
connected K\"{a}hler manifold $(N^{2n},h,J_{N})$, with connected fibers
equal to the lifted Reeb leaves. The transverse complex structures must be
compatible, and 
\begin{equation}
q^{\ast }g^{T}=D^{\ast }h,\qquad q^{\ast }f=D^{\ast }f_{N}
\label{eq:development}
\end{equation}%
for a single-valued smooth function $f_{N}$ on $N$.
\end{definition}

The connected-fiber condition ensures that basic functions descend to $N$;
local sections imply smoothness. We do not build uniqueness into the
definition. When the covering is the Riemannian universal cover, however,
the base is canonically the quotient by the lifted Reeb leaves: hence two
such developments on the universal cover are related by a unique K\"{a}hler
isometry intertwining the development maps. Accordingly, statements that do
not specify the universal covering will refer to \emph{a} complete
transverse development. The steady equation then descends by its local
transverse expression. Regular quotients and developable quasi-regular
quotients give such developments after passage to their smooth universal
covers. When imposing curvature decay, we always specify the space on which
distance is measured.

The global de Rham theorem applies on $N$, whereas a parallel transverse
splitting alone gives only local product charts. The foliated decomposition
results of \cite{Molino,BlumenthalHebda} do not by themselves supply this
complete smooth base. The required global constructions are given in
Sections~\ref{sec:quasi-regular} and~\ref{sec:irregular}.

\subsection{Sasaki-Ricci Soliton}

\begin{definition}
A Sasakian manifold $(M,\xi ,\eta ,g,\Phi )$ is $\eta $-Einstein if there is
a constant $A$ such that the Ricci curvature 
\begin{equation*}
\func{Ric}=Ag+(2n-A)\eta \otimes \eta .
\end{equation*}
\end{definition}

For $g=g^{T}+\eta\otimes\eta$, we have 
\begin{equation*}
\func{Ric}^{T}=(A+2)g^{T}
\end{equation*}
and 
\begin{equation*}
\func{Ric}=Ag^{T}+2n(\eta\otimes\eta).
\end{equation*}
Writing $\kappa=A+2$, we have 
\begin{equation*}
\rho^{T}=\kappa\,\omega^{T}.
\end{equation*}
which is also called transverse K\"{a}hler-Einstein. For $A=2n$%
\begin{equation*}
\func{Ric}=2ng^{T}+2n(\eta\otimes\eta)=2ng
\end{equation*}
which is Sasaki-Einstein.

In summary, $(M,g)$ is Sasaki-Einstein with 
\begin{equation*}
\func{Ric}_{g}=2ng
\end{equation*}%
if and only if the K\"{a}hler cone $(C(M),\overline{g})$ is Ricci-flat. On
the other hand, the transverse K\"{a}hler structure to the Reeb foliation $%
F_{\xi }$ is K\"{a}hler-Einstein with 
\begin{equation*}
\func{Ric}^{T}=2(n+1)g^{T}.
\end{equation*}%
Now we have some useful relations for curvature tensors in case of Sasakian
manifold. 
\begin{equation}
\begin{array}{lll}
\func{Rm}(e_{i},\xi ,\xi ,e_{j}) & = & g(e_{i},e_{j})=\delta _{ij}.%
\end{array}
\label{eqn-SR3}
\end{equation}%
This follows directly from the Sasakian identity %
\eqref{eq:sasaki-differential} together with $\nabla^g_X\xi=-\Phi X$; see 
\cite[Chapter~7]{BoyerGalicki}. Therefore 
\begin{equation}
\begin{array}{c}
\func{Ric}(\xi ,\xi )=\sum_{i=1}^{2n}\func{Rm}(e_{i},\xi ,\xi ,e_{i})=2n.%
\end{array}
\label{eqn-SR4}
\end{equation}%
By using (\ref{eqn-SR4}), we obtain%
\begin{equation*}
R^{T}=R+2n
\end{equation*}%
For a horizontal orthonormal frame $\{e_i\}_{i=1}^{2n}$, write
$R_{ii}=\func{Ric}_g(e_i,e_i)$, where $\func{Ric}_g$ denotes the ambient
Ricci tensor. Then
\begin{equation}
\begin{array}{c}
\sum_{i=1}^{2n}R_{ii}=\sum_{i=1}^{2n}\func{Ric}_g(e_{i},e_{i})=R^{T}-4n.%
\end{array}
\label{eqn-SR6}
\end{equation}

\begin{definition}
A gradient Sasaki--Ricci soliton is a Sasakian manifold $(M,g,\eta,\xi,\Phi)$
together with a real basic function $f$ satisfying 
\begin{equation}
\func{Ric}^{T}+\func{Hess}^{T}f=\kappa g^{T}  \label{pr1}
\end{equation}
for a constant $\kappa$. It is called expanding, steady or shrinking
according as $\kappa<0$, $\kappa=0$ or $\kappa>0$, respectively. In the
notation above, $\kappa=A+2$. The equation implies that $\func{Hess}^{T}f$
is $J^T$-invariant, or equivalently that the transverse $(1,0)$-gradient $%
\nabla^{T,1,0}f$ is transversely holomorphic. This gradient definition is
the only soliton formulation used in this paper.
\end{definition}

\begin{remark}
For a constant $a>0$, the \emph{transverse $D$-homothetic deformation} of $%
(M,\xi ,\eta ,\Phi ,g)$ is the Sasakian structure 
\begin{equation}
\eta _{a}=a\eta ,\qquad \xi _{a}=a^{-1}\xi ,\qquad \Phi _{a}=\Phi ,\qquad
g_{a}=ag+a(a-1)\eta \otimes \eta .  \label{eq:pre-Dhom}
\end{equation}%
Then $g_{a}^{T}=ag^{T}$, the transverse Levi-Civita connection is unchanged,
and $\mathrm{Ric}_{a}^{T}=\mathrm{Ric}^{T}$, $\mathrm{Hess}_{a}^{T}f=\mathrm{%
Hess}^{T}f$ as $(0,2)$-tensors. In particular, the steady equation is
preserved by a fixed $D$-homothety (with the same potential as a function),
while the Reeb period and transverse lengths are rescaled accordingly. A
fixed $D$-homothety preserves completeness and the regular, quasi-regular or
irregular type of the Reeb foliation. Indeed $g_a$ and $g$ are uniformly
equivalent by positive constants depending only on the fixed number $a$, so
they have the same Cauchy sequences and completeness is unchanged.
\end{remark}

In the remainder of the paper we consider the steady case. The potential is
a real basic function $f$, and the steady Sasaki--Ricci soliton equation is 
\begin{equation}
\mathrm{Ric}^{T}+\mathrm{Hess}^{T}f=0, \qquad (\mathrm{Hess}^{T}f)^{2,0}=(%
\mathrm{Hess}^{T}f)^{0,2}=0.  \label{eq:steady}
\end{equation}
The type condition is equivalent to transverse holomorphicity of $%
\nabla^{T,1,0}f$. With $d^{c}u=-du\circ J$, the form version is 
\begin{equation}
\rho^{T}+\tfrac12 d_{B}d_{B}^{c}f=0, \qquad \rho^{T}=\mathrm{Ric}%
^{T}(J\cdot,\cdot).  \label{eq:pre-form}
\end{equation}
The transverse trace and divergence give 
\begin{equation}
R^{T}+\Delta^{T}f=0, \qquad \nabla^{T}R^{T}=2\mathrm{Ric}^{T}(\nabla^{T}f,%
\cdot), \qquad R^{T}+|\nabla^{T}f|^{2}=\mathrm{const}.
\label{eq:pre-identities}
\end{equation}
These are the conventions used throughout the classification below.

\subsection{Soliton Identities and Curvature Estimates}

\label{subsec:soliton-identities}

\begin{lemma}
\label{lem:curvature-control} A K\"{a}hler curvature tensor of complex
dimension $n$ with nonnegative bisectional curvature satisfies $|\func{Rm}%
|\leq A_{n}R$ for a dimensional constant $A_{n}$. In particular, $R=0$
implies $\func{Rm}=0$.
\end{lemma}

\begin{proof}
Nonnegative bisectional curvature implies $\func{Ric}\geq 0$. If $R=0$, then 
$\func{Ric}=0$, and for each $(1,0)$ vector $Z$ the nonnegative Hermitian
form $\func{Rm}(Z,\overline{Z},\,\cdot \,,\overline{\,\cdot \,})$ has zero
trace. It therefore vanishes, and polarization gives $\func{Rm}=0$. The
scalar trace consequently has a positive minimum on the compact set of
algebraic K\"{a}hler curvature tensors with nonnegative bisectional
curvature and unit norm. The estimate follows by rescaling.
\end{proof}

\begin{theorem}
\label{thm:conservation} On a connected gradient steady Sasaki-Ricci soliton
there is a real constant $C_{0}$ such that 
\begin{equation}
R^{T}+|\nabla ^{T}f|^{2}=C_{0},  \label{eq:conservation}
\end{equation}%
and the drift operator $L_{f}$ satisfies 
\begin{equation}
L_{f}R^{T}=-2|\func{Ric}^{T}|^{2}.  \label{eq:drift}
\end{equation}%
These identities require neither completeness nor a regular Reeb flow. If $%
B^{T}\geq 0$, then 
\begin{equation}
0\leq R^{T}\leq C_{0},\qquad |\nabla ^{T}f|\leq \sqrt{C_{0}},
\label{eq:bounded-gradient}
\end{equation}%
and either $R^{T}>0$ everywhere or $\func{Rm}^{T}\equiv 0$. In the latter
case, if $g$ is complete, its Riemannian universal cover is the standard
Heisenberg Sasakian manifold.
\end{theorem}

\begin{proof}
The usual steady soliton identities apply in each transverse K\"{a}hler
chart, see \cite[Section~4]{DengZhu2015} (with potential $-f$) and \cite[%
Section~2, especially Lemma~2.1]{PetersenWylie}. In our conventions these
yield \eqref{eq:drift} together with 
\begin{equation}
\Delta ^{T}f=-R^{T},\qquad \nabla ^{T}R^{T}=2\func{Ric}^{T}(\nabla
^{T}f,\,\cdot \,).  \label{eq:gradient-scalar}
\end{equation}%
Differentiating $R^{T}+|\nabla ^{T}f|^{2}$ and using the soliton equation
shows that its transverse differential vanishes. Since the function is basic
and $M$ is connected, it is constant. The bounds follow from $R^{T}\geq 0$.
By \eqref{eq:basic-ambient}, the strong minimum principle for the ambient
drift operator applies to $R^{T}$ and gives the stated dichotomy. Lemma~\ref%
{lem:curvature-control} yields transverse flatness when $R^{T}\equiv 0$. For
every unit horizontal vector $X$, the Sasakian curvature identity 
\begin{equation*}
K^T(X,J^TX)=K_g(X,\Phi X)+3
\end{equation*}
(see \cite[Chapter~7]{BoyerGalicki}) shows that transverse flatness implies $%
K_g(X,\Phi X)=-3$. Hence the full Sasakian metric has constant $\Phi$%
-holomorphic sectional curvature $-3$, and \cite[Proposition~4.1(ii)]%
{Tanno1969} identifies its complete simply connected cover with the
Heisenberg model.
\end{proof}

\begin{remark}
A transverse critical point gives $C_{0}=\max_{M}R^{T}$. Transverse
flatness, however, allows affine potentials on the Heisenberg cover,
possibly with $C_{0}>\sup R^{T}$.
\end{remark}

\begin{lemma}
\label{lem:scalar-gradient} A gradient steady soliton with $\func{Ric}\geq 0$
satisfies 
\begin{equation}
|\nabla R|^{2}\leq 4C_{0}R^{2}.  \label{eq:scalar-gradient-bound}
\end{equation}%
The same estimate holds for transverse tensors on a Sasakian soliton.
\end{lemma}

\begin{proof}
Use $\nabla R=2\func{Ric}(\nabla f,\,\cdot \,)$. Since $\func{Ric}\geq0$,
its operator norm is its largest eigenvalue and is at most its trace $R$;
hence $\Vert\func{Ric}\Vert_{\mathrm{op}}\leq R$. Together with $|\nabla
f|^{2}\leq C_{0}$ this gives the claim.
\end{proof}

\section{Regular Solitons and Complete Transverse Models}

\label{sec:regular}

\subsection{Complete Quotients and Splitting}

\label{subsec:regular-splitting}

For a regular Reeb circle action or a proper Reeb line action, the quotient 
\begin{equation*}
\pi :M\rightarrow Q
\end{equation*}
is a complete K\"{a}hler manifold and the potential satisfies %
\eqref{eq:quotient-soliton}. We first describe the connection data of the
lift. Passing to the universal cover of $Q$ then reduces the transverse
classification to the simply connected K\"{a}hler case.

\begin{proposition}
\label{prop:regular-lift-data} Let $(Q,h,J,\omega ,f_{Q})$ be a complete K%
\"{a}hler steady soliton.

\begin{enumerate}
\item A regular Sasakian circle lift with Reeb period $L>0$ exists precisely
when $2[\omega ]/L$ is the real image of a class in $H^{2}(Q;\mathbb{Z})$.
It is specified by a principal $\mathbb{R}/L\mathbb{Z}$-bundle and a
connection $\eta $ with $d\eta =2\pi ^{\ast }\omega $. Its metric and
potential are 
\begin{equation*}
g=\pi ^{\ast }h+\eta ^{2}\text{ \ \ \ \textrm{and} \ \ }f=\pi ^{\ast }f_{Q}.
\end{equation*}

\item Or a proper-line lift exists precisely when $\omega $ is exact. In
that case it is given on 
\begin{equation*}
Q\times \mathbb{R}
\end{equation*}
with 
\begin{equation*}
\eta =dt+A\text{ \ \ \textrm{and} \ }dA=2\omega .
\end{equation*}%
All these lifts are complete, and their transverse curvatures equal those of 
$h$.
\end{enumerate}
\end{proposition}

\begin{proof}
Write $d\eta =\pi ^{\ast }\Omega $. For a circle connection of period $L$,
Chern--Weil theory gives 
\begin{equation*}
e_{\mathbb{R}}=[\Omega /L]=2[\omega ]/L.
\end{equation*}%
Conversely, choose a bundle with this integral class and a connection $\eta
_{0}$ of curvature $\pi ^{\ast }\Omega $. Since $[\Omega ]=2[\omega ]$,
there is a one-form $a$ on $Q$ with $\Omega -2\omega =da$. Then $\eta =\eta
_{0}-\pi ^{\ast }a$ has the required curvature. Every principal $\mathbb{R}$%
-bundle is trivial, so a line connection has the form $dt+A$; its curvature
is $2\omega $ precisely when $dA=2\omega $.

For either lift, put $\Phi =-J^{\mathrm{hor}}$ on $\ker \eta $ and $\Phi \xi
=0$. Koszul's formula gives, for horizontal lifts, 
\begin{equation*}
\nabla _{X}^{g}Y=(\nabla _{X}^{h}Y)^{\mathrm{hor}}-\omega (X,Y)\xi ,\qquad
\nabla _{X}^{g}\xi =JX.
\end{equation*}%
Together with $\nabla ^{h}J=0$, these imply the Sasakian identity $(\nabla
_{X}^{g}\Phi )Y=g(X,Y)\xi -\eta (Y)X$. The transverse metric is $h$, so the
steady equation and curvature assertions follow by pullback.

To prove completeness, let $c=(x,t)$ be a curve of finite length. Its
projection converges in $Q$. In a local trivialization near the limit, the
connection coefficients are bounded and 
\begin{equation*}
|\dot{t}|\leq |\eta (\dot{c})|+|A(\dot{x})|\leq (1+\sup |A|_{h})|\dot{c}%
|_{g}.
\end{equation*}%
Thus the fiber coordinate also converges. For a circle fiber the same
estimate is applied to a lift of its coordinate to $\mathbb{R}$. Hence no
finite-length curve escapes to infinity.
\end{proof}

\begin{theorem}
\label{thm:kahler-splitting} Let $(Q^{2n},h,J,f)$ be a complete simply
connected gradient steady K\"{a}hler-Ricci soliton with 
\begin{equation*}
B_{h}\geq 0
\end{equation*}%
satisfying 
\begin{equation*}
\func{Ric}_{h}+\func{Hess}_{h}f=0.
\end{equation*}%
Then there is a holomorphic isometric decomposition 
\begin{equation}
(Q,h,J)\cong (P_{1}^{2k_{1}},h_{1},J_{1})\times \cdots \times
(P_{s}^{2k_{s}},h_{s},J_{s})\times (\mathbb{C}^{m},g_{\mathrm{flat}%
},J_{0}),\qquad \sum_{i=1}^{s}k_{i}+m=n.  \label{eq:complete-model-splitting}
\end{equation}%
Each $P_{i}$ is complete, simply connected, noncompact and irreducible, with
nonnegative holomorphic bisectional curvature and positive Ricci curvature.
The potential has the form 
\begin{equation}
f=\sum_{i=1}^{s}f_{i}\circ \pi _{i}+\ell \circ \pi _{\mathrm{flat}}+c,
\label{eq:potential-product}
\end{equation}%
where $\func{Ric}_{h_{i}}+\func{Hess}_{h_{i}}f_{i}=0$ and $\ell $ is a real
linear function on $\mathbb{C}^{m}\cong \mathbb{R}^{2m}$. Empty products and 
$m=0$ are allowed. If $f$ has a critical point, then $\ell =0$.
\end{theorem}

\begin{proof}
The conservation identity and Lemma~\ref{lem:curvature-control} give 
\begin{equation}
0\leq R_{h}\leq C_{0},\qquad |\nabla f|\leq \sqrt{C_{0}},\qquad |\func{Rm}%
_{h}|\leq c_{n}R_{h}\leq c_{n}C_{0}.  \label{eq:model-bounded-curvature}
\end{equation}%
The flat case is Euclidean with affine potential. In the nonflat case,
integration of $\Delta f=-R_{h}$ rules out compactness. The bounded
holomorphic field $\nabla f$ has a complete flow $\psi _{t}$, and 
\begin{equation*}
h(t)=\psi _{t}^{\ast }h,\qquad \partial _{t}h(t)=-2\func{Ric}_{h(t)},\qquad
-\infty <t<\infty .
\end{equation*}%
Cao's theorem \cite[Theorem~2.1]{Cao2004} splits off a flat factor whose
tangent bundle is the Ricci kernel. The complementary factor has positive
Ricci curvature. Its de Rham decomposition gives %
\eqref{eq:complete-model-splitting}; since $J$ is parallel, each nonflat
factor is K\"{a}hler and the splitting is holomorphic \cite{Wu}.

It remains to decompose the potential. For vector fields $X_i,X_j$ lifted
from distinct factors, 
\begin{equation*}
X_iX_jf=\func{Hess} f(X_i,X_j)=-\func{Ric}(X_i,X_j)=0.
\end{equation*}
Hence $f$ is a sum of factor potentials. Its flat summand has zero Hessian
and is affine, giving \eqref{eq:potential-product}. A compact nonflat factor
is impossible, since 
\begin{equation*}
0=\int_{P_i}\Delta f_i\,dV=-\int_{P_i}R_{h_i}\,dV<0.
\end{equation*}
Finally, at a critical point of $f$ every factor component of its gradient
vanishes. In particular, the constant flat gradient is zero, so $\ell=0$.
\end{proof}

\subsection{The Five-Dimensional Transverse Structure}

\label{subsec:five-dimensional-structure}

\begin{theorem}
\label{thm:five} Let $(M^{5},g,\eta ,\xi ,\Phi ,f)$ be a complete gradient
steady Sasaki--Ricci soliton with 
\begin{equation*}
B^{T}\geq 0.
\end{equation*}
Assume that its Reeb foliation and potential admit a complete smooth
transverse development $(N^{4},h,J,f_{N})$ in the sense of Definition~\ref%
{def:development}. Then $N$ is holomorphically isometric to precisely one of
the following types:

\begin{enumerate}
\item the flat space $(\mathbb{C}^{2},g_{\mathrm{flat}})$,

\item $(\mathbb{C},ag_{\mathrm{cig}})\times (\mathbb{C},g_{\mathrm{flat}})$,
for some $a>0$,

\item $(\mathbb{C},ag_{\mathrm{cig}})\times (\mathbb{C},bg_{\mathrm{cig}})$,
for some $a,b>0$,

\item a complete simply connected irreducible K\"{a}hler surface with
nonnegative bisectional curvature and positive Ricci curvature, carrying a
gradient steady soliton.
\end{enumerate}

Here $g_{\mathrm{cig}}=|\mathrm{d}z|^{2}/(1+|z|^{2})$. Case \textup{(4)}
includes Cao's $U(2)$-invariant metric and the non-Cao examples of
Proposition~\ref{prop:ac}.
\end{theorem}

\begin{lemma}
\label{lem:surface} A complete simply connected nonflat gradient steady
Ricci soliton surface with nonnegative Gaussian curvature is a positive
multiple of the cigar.
\end{lemma}

\begin{proof}
The scalar minimum principle gives $K>0$. For $F=-f$, 
\begin{equation*}
\func{Hess}F=Kh,\qquad dK=-K\,dF,\qquad 2K+|\nabla F|^{2}=C_{0}>0.
\end{equation*}%
If $F$ has no critical point, along a complete trajectory of $-\nabla F$ one
obtains 
\begin{equation*}
K^{\prime }=K(C_{0}-2K),\qquad K(t)=\frac{C_{0}}{2(1+be^{-C_{0}t})},\qquad
|\nabla F|^{2}=\frac{C_{0}be^{-C_{0}t}}{1+be^{-C_{0}t}},\quad b>0.
\end{equation*}%
The speed is integrable on $[0,\infty )$, so completeness gives a limit at
which $\nabla F=0$, a contradiction. Strict convexity makes this minimum $o$
unique. Normalize $F(o)=0$ and choose $a,\mu >0$ with $\func{Hess}F\geq \mu
h $ on $\overline{B}(o,a)$. Along minimizing geodesics from $o$, 
\begin{equation*}
F(x)\geq \mu a\bigl(d_{h}(o,x)-a\bigr)\qquad \bigl(d_{h}(o,x)\geq a\bigr).
\end{equation*}%
Thus $F$ is proper.

The field $J\nabla F$ is Killing, and its linearization at $o$ is rotation
with speed $c=K(o)>0$. Its time-$2\pi /c$ map fixes $o$ with identity
differential and is therefore the identity isometry. The positive levels of $%
F$ are connected circles: this holds near $o$ by the Morse lemma and extends
to all levels by the flow of $\nabla F/|\nabla F|^{2}$. Hence these levels
are the Killing orbits. Moreover, $\nabla F/|\nabla F|$ is geodesic, giving
global polar coordinates 
\begin{equation*}
h=ds^{2}+a(s)^{2}d\theta ^{2},\qquad a(0)=0,\quad a^{\prime }(0)=1.
\end{equation*}%
The soliton equation becomes 
\begin{equation*}
F^{\prime \prime }=\frac{F^{\prime }a^{\prime }}{a}=-\frac{a^{\prime \prime }%
}{a},\qquad F^{\prime }=ca,\qquad a^{\prime \prime }+caa^{\prime }=0.
\end{equation*}%
Since 
\begin{equation*}
\frac{d}{ds}\left( a^{\prime }+\frac{c}{2}a^{2}\right) =a^{\prime \prime
}+caa^{\prime }=0,
\end{equation*}%
there is a constant $C$. Smooth polar coordinates give $a(0)=0$ and $%
a^{\prime }(0)=1$, hence $C=1$ and 
\begin{equation}
a^{\prime }+\frac{c}{2}a^{2}=1.  \label{eq:cigar-first-integral}
\end{equation}%
This Riccati equation has the unique solution with $a(0)=0$, 
\begin{equation*}
a(s)=\sqrt{\frac{2}{c}}\tanh \!\left( \sqrt{\frac{c}{2}}s\right) .
\end{equation*}%
Indeed direct differentiation gives $a^{\prime }=\func{sech}^{2}(\sqrt{c/2}%
\,s)$, so $a^{\prime }(s)+\frac{c}{2}a(s)^{2}=1$ and therefore $a^{\prime
\prime }+caa^{\prime }=0$. At the tip $|\nabla F|=0$ and $K(o)=c$, so the
conservation law $2K+|\nabla F|^{2}=C_{0}$ gives $C_{0}=2c.$ No additional
normalization of $c$ is being imposed. With $r=\sinh (\sqrt{c/2}\,s)$ this
gives $h=(2/c)(dr^{2}+r^{2}d\theta ^{2})/(1+r^{2})$, the corresponding
scaling of Hamilton's cigar (cf.\ \cite{Hamilton1988}).
\end{proof}

\begin{proof}[Proof of Theorem~\protect\ref{thm:five}]
Apply Theorem~\ref{thm:kahler-splitting} in complex dimension two. The
complex dimensions of the nonflat factors are respectively $0$, $1$, or $2$.
Lemma~\ref{lem:surface} identifies every one-dimensional nonflat factor with
a scaled cigar.
\end{proof}

\subsection{Potential Levels and Curvature Estimates}

\label{subsec:potential-levels}

Let $(N^{2k},h,J,f_N)$ be a complete K\"ahler steady soliton with positive
Ricci curvature. Once a critical point $o$ is fixed, we use the increasing
potential 
\begin{equation}  \label{eq:F-potential}
F=f_N(o)-f_N,\qquad \func{Hess} F=\func{Ric}_h.
\end{equation}

\begin{proposition}
\label{prop:potential} Let $(N,h,J,f_{N})$ be a complete noncompact gradient
steady K\"{a}hler-Ricci soliton with 
\begin{equation*}
\func{Ric}_{h}>0.
\end{equation*}%
Suppose%
\begin{equation*}
R_{h}(x)\rightarrow 0\text{ \ \ \textrm{as} \ \ \ }d_{h}(o_{\ast
},x)\rightarrow \infty .
\end{equation*}%
Then

\begin{enumerate}
\item $f_{N}$ has a unique critical point $o$, which is its maximum, and 
\begin{equation*}
C_{0}=R_{h}(o)>0,\qquad F\geq 0,\qquad F(o)=0.
\end{equation*}

\item There are constants $c,C>0$, depending only on the fixed soliton and
the choice of the compact ball used below, such that 
\begin{equation}
c\rho (x)-C\leq F(x)\leq \sqrt{C_{0}}\rho (x),\qquad \rho (x)=d_{h}(o,x).
\label{eq:potential-linear}
\end{equation}

\item In particular $F$ is proper. Every level $\Sigma _{s}=\{F=s\}$, $s>0$,
is a smooth compact hypersurface diffeomorphic to $S^{2k-1}$.
\end{enumerate}
\end{proposition}

\begin{proof}
Since $R_{h}>0$ and $R_{h}$ tends to zero at infinity, it attains a positive
maximum. At a maximum point $o$, 
\begin{equation*}
0=\nabla R_{h}(o)=2\func{Ric}_{h}(\nabla f_{N},\cdot )(o),\qquad \nabla
f_{N}(o)=0,\qquad R_{h}(o)=C_{0}.
\end{equation*}%
For a unit-speed minimizing geodesic $\gamma $ from $o$, 
\begin{equation*}
(F\circ \gamma )^{\prime }(0)=0,\qquad (F\circ \gamma )^{\prime \prime }=%
\func{Ric}_{h}(\dot{\gamma},\dot{\gamma})>0.
\end{equation*}%
Thus $o$ is the unique critical point and the absolute minimum of $F$.
Choose $a,\mu >0$ with $\func{Ric}_{h}\geq \mu h$ on $\overline{B}(o,a)$.
For $\rho =d(o,x)\geq a$, integration gives 
\begin{equation*}
\mu a(\rho -a)\leq \int_{0}^{\rho }(\rho -t)\func{Ric}_{h}(\dot{\gamma},\dot{%
\gamma})\,dt=F(x)\leq \sqrt{C_{0}}\rho .
\end{equation*}%
This proves \eqref{eq:potential-linear}; Hopf--Rinow then gives properness.
At $o$ the Hessian is positive definite, so the Morse lemma identifies small
positive levels with $S^{2k-1}$. On each compact slab between positive
levels, the field $V=\nabla F/|\nabla F|^{2}$ is smooth and satisfies $VF=1$%
. Its flow identifies all positive levels with these small spheres.
\end{proof}

Once a critical point exists, the growth and level-set arguments require
only $\func{Ric}_{h}>0$. Scalar decay is used above to produce that point.

\begin{remark}
On a Sasakian lift $f=\pi ^{\ast }f_{N}$, the fiber over $o$ is transversely
critical. For a line contactization, $\{f=f_{N}(o)-s\}=\Sigma _{s}\times 
\mathbb{R}$ is noncompact. Compactness of $\ \Sigma _{s}$ concerns the base.
\end{remark}

Under the stronger estimate $R_h\leq C/\rho$, the comparison %
\eqref{eq:potential-linear} implies, for large $s$, 
\begin{equation}  \label{eq:level-bounds}
\rho\asymp s\text{ on }\Sigma_s,\qquad R_h\leq\frac{C_1}{s},\qquad |\nabla
F|^2=C_0-R_h\geq\frac{C_0}{2}.
\end{equation}
We write $d=2k$, $h_s=h|_{T\Sigma_s}$ and $\Omega_s=\{F\leq s\}$. The next
proposition records the geometry of these level hypersurfaces.

\begin{proposition}
\label{prop:level-geometry} With outward normal $\nu=\nabla F/|\nabla F|$,
the second fundamental form of $\partial\Omega_s=\Sigma_s$, defined by $%
A_s(X,Y)=h(\nabla_X\nu,Y)$, satisfies 
\begin{equation}  \label{eq:level-second-form}
A_s(X,Y)=\frac{\func{Ric}_h(X,Y)}{|\nabla F|}.
\end{equation}
Let $V=\nabla F/|\nabla F|^2$. Its flow identifies the positive levels of $F$%
, and on tangent vectors to a level, 
\begin{equation}  \label{eq:level-metric-evolution}
(\mathcal{L}_Vh)(X,Y) =\frac{2\func{Ric}_h(X,Y)}{|\nabla F|^2}.
\end{equation}
If $R_h(x)\leq C/\rho(x)$ for large $\rho(x)=d_h(o,x)$, then for all
sufficiently large $s$, 
\begin{equation}  \label{eq:level-shape-bound}
R_h\leq\frac{C^{\prime }}s,\qquad |\nabla F|^2\geq\frac{C_0}{2},\qquad
|A_s|\leq\frac{C^{\prime \prime }}s \quad\hbox{on }\Sigma_s.
\end{equation}
\end{proposition}

\begin{proof}
For $X,Y\in T\Sigma_s$ and $a=|\nabla F|^{-2}$, 
\begin{equation*}
\begin{aligned} h(\nabla_X\nu,Y)&=|\nabla F|^{-1}\operatorname{Hess}
F(X,Y),\\ \mathcal L_{a\nabla F}h &=2a\operatorname{Hess} F+da\otimes
dF+dF\otimes da. \end{aligned}
\end{equation*}
Since $dF(X)=dF(Y)=0$ and $\func{Hess} F=\func{Ric}_h$, these prove %
\eqref{eq:level-second-form}--\eqref{eq:level-metric-evolution}. The flow
exists across compact positive slabs and $VF=1$. Finally, %
\eqref{eq:level-bounds} and $|\func{Ric}_h|\leq R_h$ give %
\eqref{eq:level-shape-bound}.
\end{proof}

For $h_s$-orthonormal tangent vectors $X,Y$, the Gauss equation is 
\begin{equation}  \label{eq:level-gauss}
K_{h_s}(X,Y)=K_h(X,Y) +A_s(X,X)A_s(Y,Y)-A_s(X,Y)^2.
\end{equation}
Consequently, if $|\func{Rm}_h|\leq C/\rho$ outside a compact set, then $%
|K_{h_s}|\leq C^{\prime }/s$. Rescaling by $s^{-1}$ multiplies sectional
curvature by $s$, and therefore 
\begin{equation}  \label{eq:scaled-level-curvature}
\sup_{\Sigma_s}|K_{s^{-1}h_s}|\leq C^{\prime \prime }.
\end{equation}

For the rest of this subsection we impose 
\begin{equation}  \label{eq:level-dz-hypotheses}
\func{Sec}_h\geq0,\qquad \func{Ric}_h>0, \qquad R_h(x)\leq\frac C{\rho(x)}
\quad(\rho(x)\geq\rho_0).
\end{equation}
In particular, $|\func{Rm}_h|\leq c_dR_h\leq C^{\prime }/\rho$. We use the
following estimates of Deng--Zhu.

\begin{proposition}
\label{prop:level-diameter} Under \eqref{eq:level-dz-hypotheses}, there are $%
D,s_{0}>0$ such that 
\begin{equation}
\func{diam}(\Sigma _{s},h_{s})\leq D\sqrt{s},\qquad \func{diam}(\Sigma
_{s},s^{-1}h_{s})\leq D\quad (s\geq s_{0}).  \label{eq:level-diameter}
\end{equation}
\end{proposition}

\begin{proof}
Identify the levels by the flow of $V$, and denote their intrinsic distances
by $d_{s}$. With the convention $\func{Hess}F=\func{Ric}_{h}$, \cite[%
Lemmas~2.4, 3.2 and Proposition~3.3]{DengZhu2020} apply under %
\eqref{eq:level-dz-hypotheses}. Their second-variation estimate gives 
\begin{equation*}
\frac{d^{+}}{ds}d_{s}(x,y)\leq C\left( \frac{\tau }{s}+\frac{1}{\tau }+\frac{%
d_{s}(x,y)}{s^{2}}\right) ,\qquad d_{s}(x,y)\geq 2\tau .
\end{equation*}%
Taking $\tau =\sqrt{s}$ yields $D^{+}d_{s}\leq Cs^{-1/2}+Cs^{-2}d_{s}$
whenever $d_{s}\geq 2\sqrt{s}$. Integration gives $d_{s}\leq D\sqrt{s}$, and
distances below the threshold already have this bound. Since diameters scale
by $s^{-1/2}$ under $h_{s}\mapsto s^{-1}h_{s}$, the second assertion follows.
\end{proof}

\begin{lemma}
\label{lem:vanishing-propagation} Assume \eqref{eq:level-dz-hypotheses}. If $%
p_{i}\in \Sigma _{s_{i}}$ satisfy $s_{i}\rightarrow \infty $ and $%
s_{i}R_{h}(p_{i})\rightarrow 0$, then 
\begin{equation}
s_{i}\sup_{\Sigma _{s_{i}}}R_{h}\longrightarrow 0.
\label{eq:uniform-level-vanishing}
\end{equation}
\end{lemma}

\begin{proof}
By \cite[Lemma~4.4]{DengZhu2020}, there is a fixed $\varepsilon >0$ such
that 
\begin{equation}
\sup_{B_{h}(p_{i},\varepsilon \sqrt{s_{i}})}s_{i}R_{h}\longrightarrow 0.
\label{eq:local-zero-propagation}
\end{equation}%
Their formulation with $FR_{h}$ is equivalent, since on this ball $%
|F-s_{i}|\leq \sqrt{C_{0}}\varepsilon \sqrt{s_{i}}$. For any $q_{i}\in
\Sigma _{s_{i}}$, \eqref{eq:level-diameter} joins $p_{i}$ to $q_{i}$ by an
intrinsic geodesic of length at most $D\sqrt{s_{i}}$. Divide it into at most 
$\lceil 2D/\varepsilon \rceil +1$ segments and apply %
\eqref{eq:local-zero-propagation} successively at their endpoints. Then $%
s_{i}R_{h}(q_{i})\rightarrow 0$. Taking $q_{i}$ at the scalar maximum proves %
\eqref{eq:uniform-level-vanishing}; cf. \cite[Lemma~4.5]{DengZhu2020}.
\end{proof}

\begin{theorem}
\label{thm:factor-decay} Let $d\geq 3$ and assume %
\eqref{eq:level-dz-hypotheses}. Then there exist $c,C_{1},\rho _{1}>0$ such
that 
\begin{equation}
\frac{c}{\rho (x)}\leq R_{h}(x)\leq \frac{C_{1}}{\rho (x)}\qquad (\rho
(x)\geq \rho _{1}).  \label{eq:factor-two-sided}
\end{equation}
\end{theorem}

\begin{proof}
This is \cite[Proposition~4.3]{DengZhu2020}. If the lower bound failed,
choose $p_{i}\rightarrow \infty $ with $\rho (p_{i})R_{h}(p_{i})\rightarrow 0
$ and put $s_{i}=F(p_{i})$. Since $F\asymp \rho $, the preceding lemma gives $%
\varepsilon _{i}=s_{i}\sup_{\Sigma _{s_{i}}}R_{h}\rightarrow 0$. The Gauss
equation, $|\func{Rm}_{h}|\leq c_{d}R_{h}$, and \eqref{eq:level-second-form}
imply 
\begin{equation*}
\sup_{\Sigma _{s_{i}}}|K_{s_{i}^{-1}h_{s_{i}}}|\leq C(\varepsilon
_{i}+\varepsilon _{i}^{2}/s_{i})\longrightarrow 0.
\end{equation*}%
Together with \eqref{eq:level-diameter}, this gives 
\begin{equation*}
\sup |K_{s_{i}^{-1}h_{s_{i}}}|\,\func{diam}(\Sigma
_{s_{i}},s_{i}^{-1}h_{s_{i}})^{2}\longrightarrow 0.
\end{equation*}%
Gromov's almost-flat theorem \cite{Gromov} implies that the universal cover
of each sufficiently distant level is diffeomorphic to $\mathbb{R}^{d-1}$.
This contradicts $\Sigma _{s_{i}}\cong S^{d-1}$, since $d-1\geq 2$. The
upper bound is part of the hypothesis.
\end{proof}

In real dimension two the level is a circle, whose universal cover is $%
\mathbb{R}$. Thus the topological contradiction does not apply. This is
consistent with the exponential scalar curvature decay of the cigar. For $%
d\geq 3$, \eqref{eq:factor-two-sided} and $F\asymp \rho $ also give a uniform
comparison of the maximum and minimum of $R_{h}$ on every sufficiently
distant level.

\subsection{Curvature Decay on Quotients and Developments}

\label{subsec:quotient-decay}

We first record the classification under a decay bound on a smooth quotient.
We then apply it to the quotient of the universal cover by the lifted Reeb
flow, using the decay bound of Theorem~\ref{thm:intro-decay}.

\begin{theorem}
\label{thm:decay} Let a complete gradient steady Sasaki-Ricci soliton have a
smooth complete K\"{a}hler quotient $(Q^{2n},h,J,f_{Q})$ as in Section~\ref%
{sec:geometry}. Assume 
\begin{equation}
\func{Sec}_{h}\geq 0,\qquad R_{h}(x)\leq \frac{C}{\rho _{Q}(x)}\quad
\label{eq:quotient-decay}
\end{equation}%
for $\rho _{Q}(x)=d_{h}(o,x)\geq r_{0}.$ Then the universal cover of $Q$ is
holomorphically isometric to one of :

\begin{enumerate}
\item $(\mathbb{C}^{n},g_{\mathrm{flat}})$;

\item $(\mathbb{C},ag_{\mathrm{cig}})\times (\mathbb{C}^{n-1},g_{\mathrm{flat%
}})$, with $a>0$;

\item $(P^{2k},h_{P},J_{P})\times (\mathbb{C}^{n-k},g_{\mathrm{flat}})$,
where $2\leq k\leq n$, $P$ is complete and simply connected, 
\begin{equation*}
\func{Sec}_{h_{P}}\geq 0,\func{Ric}_{h_{P}}>0,
\end{equation*}%
and 
\begin{equation}
\frac{c}{\rho _{P}(x)}\leq R_{h_{P}}(x)\leq \frac{C^{\prime }}{\rho _{P}(x)}%
\qquad (\rho _{P}(x)\geq r_{1}),  \label{eq:factor-linear}
\end{equation}%
for positive constants $c,C^{\prime },r_{1}$ and intrinsic distance $\rho
_{P}$ on $P$. In this case $P$ is irreducible.
\end{enumerate}
\end{theorem}

\begin{proof}
If $Q$ is compact, then $0=\int_{Q}\Delta f_{Q}=-\int_{Q}R_{h}$. Nonnegative
sectional curvature therefore gives flatness. In complex dimension one, the
scalar minimum principle and Lemma~\ref{lem:surface} give the flat or cigar
alternatives.

We may therefore assume that $Q$ is noncompact and $n\geq 2$. Apply \cite[%
Theorem~1.1]{DengZhu2020} to $Q$, using the distance in %
\eqref{eq:quotient-decay}. Its universal cover is flat, a cigar times
Euclidean space, or $P^{d}\times \mathbb{R}^{2n-d}$ with $d>2$, $\func{Ric}%
_{h_{P}}>0$, $\func{Sec}_{h_{P}}\geq 0$ and two-sided linear scalar decay in
the intrinsic distance on $P$.

We verify that these Riemannian splittings are holomorphic. The lifted
complex structure is parallel, and the Ricci kernel is $J$-invariant. In
each nonflat case this kernel is exactly the Euclidean distribution. Thus $J$
preserves both factors, their real dimensions are even, and $d=2k$ in the
third case. In the flat case, a constant orthogonal complex structure
identifies the cover with standard $\mathbb{C}^{n}$.

If $P=P_1\times P_2$ were reducible, the mixed Hessian equation would split
its potential. Each factor has positive Ricci curvature and is noncompact,
since a compact steady factor satisfies $\int R=-\int\Delta f=0$. For fixed $%
p\in P_1$ and $y_j\to\infty$ in $P_2$, 
\begin{equation*}
R_{h_P}(p,y_j)=R_{h_1}(p)+R_{h_2}(y_j)\geq R_{h_1}(p)>0,
\end{equation*}
contradicting decay on $P$. Hence $P$ is irreducible.
\end{proof}

For a covering $\pi :\widetilde{Q}\rightarrow Q$, 
\begin{equation}
d_{h}(\pi \widetilde{q}_{0},\pi \widetilde{q})\leq d_{\widetilde{h}}(%
\widetilde{q}_{0},\widetilde{q}).  \label{eq:cover-distance-direction}
\end{equation}%
Thus decay in quotient distance need not lift to decay in covering distance.

\begin{proposition}
\label{prop:product-obstruction} Suppose $P_{1}\times P_{2}$ is a complete
Riemannian product, $R_{P_{1}}(p)>0$ at some point $p$, $R_{P_{2}}\geq 0$,
and $P_{2}$ is noncompact. Then the scalar curvature of the product does not
tend uniformly to zero with distance from a fixed point.
\end{proposition}

\begin{proof}
Choose $y_{j}\in P_{2}$ with $d_{P_{2}}(y_{0},y_{j})\rightarrow \infty $.
Product distance from $(p,y_{0})$ to $(p,y_{j})$ tends to infinity, whereas 
\begin{equation*}
R_{P_{1}\times P_{2}}(p,y_{j})=R_{P_{1}}(p)+R_{P_{2}}(y_{j})\geq
R_{P_{1}}(p)>0.
\end{equation*}%
Changing the fixed reference point does not affect the conclusion.
\end{proof}

Applied to Theorem~\ref{thm:kahler-splitting}, this obstruction shows that a
nonflat complete simply connected model with uniform scalar decay has only
one factor, which is irreducible and has positive Ricci curvature.

\begin{proof}[Proof of Theorem~\protect\ref{thm:intro-decay}]
The K\"{a}hler curvature symmetries and the first Bianchi identity give
\begin{equation*}
\func{Rm}^{T}(X,J^{T}X,J^{T}Y,Y)=\func{Rm}^{T}(X,Y,Y,X)
 +\func{Rm}^{T}(X,J^{T}Y,J^{T}Y,X)
\end{equation*}
for horizontal vectors $X,Y$. Both terms on the right are nonnegative by
$\func{Sec}^{T}\geq 0$. Thus $B^T\geq 0$, and
the properness theorem for the lifted Reeb flow
\cite[Theorem~1.1(2)]{Rigidity2026} gives a free proper lifted Reeb action and
a complete smooth K\"{a}hler quotient $D:\widetilde M\to N$. The fibers are
Reeb lines and $D$ is a principal $\mathbb{R}$-bundle. In particular $N$ is
simply connected. In the local transverse charts for $D$, the transverse
metric, curvature and Hessian are pulled back from $N$. The lifted basic
potential descends to $f_N$, so
\begin{equation*}
\widetilde R^T=D^*R_h,\qquad \func{Sec}_h\geq 0,\qquad
\func{Ric}_h+\func{Hess}_h f_N=0.
\end{equation*}
The bundle is trivial because $\mathbb{R}$ is contractible. In a
trivialization with $\widetilde\xi=\partial_t$, the form
$\widetilde\eta-dt$ is basic and hence equals $D^*\beta$. The identities
$d\widetilde\eta=2D^*\omega_h$ and
$\widetilde g=D^*h+\widetilde\eta^{\,2}$ give the asserted contactization.

Put $o=D(\widetilde o)$. The Riemannian submersion $D$ decreases lengths,
so $d_h(D(\widetilde x),o)\leq \rho(\widetilde x)$. For the reverse
inequality, horizontally lift a minimizing geodesic from
$D(\widetilde x)$ to $o$, starting at $\widetilde x$. The lift has the same
length, and completeness of $\widetilde g$ ensures that it exists up to the
endpoint. That endpoint lies in $D^{-1}(o)=\widetilde{\mathcal L}_{\widetilde o}$.
Thus
\begin{equation}
\rho(\widetilde x)=d_h(D(\widetilde x),o).  \label{eq:lifted-leaf-distance}
\end{equation}
It follows that \eqref{eq:lifted-leaf-decay} is equivalent to
$R_h(y)\leq C/d_h(o,y)$ for $d_h(o,y)\geq r_0$. By the triangle inequality
the choice of $o$ affects only the constants in this estimate.

Apply Theorem~\ref{thm:decay} to $\widetilde M$ with smooth quotient $N$.
Since $N$ is simply connected, its universal cover is $N$ itself. In either
nonflat alternative of that theorem the positive-Ricci factor has positive
scalar curvature. Proposition~\ref{prop:product-obstruction} excludes a
positive-dimensional Euclidean factor, since $R_h$ tends uniformly to zero
at infinity. The cigar therefore occurs only for $n=1$. In the remaining
nonflat case $P=N$ is irreducible and \eqref{eq:factor-linear} gives the
two-sided estimate on $N$.

If $N$ is flat, the soliton equation gives $\func{Hess}_h f_N=0$, so $f_N$
is affine. In the nonflat cases, $R_h>0$ and the upper decay bound and
completeness imply that $R_h$ attains a positive maximum. By
\eqref{eq:gradient-scalar}, at this point $\nabla f_N=0$ since
$\func{Ric}_h>0$. The identity $\func{Hess}_h f_N=-\func{Ric}_h$ gives
strict concavity along every nonconstant geodesic, so this point is the
unique critical point and the global maximum of $f_N$.
\end{proof}

For example, the product of the cigar with a flat complex torus satisfies
\eqref{eq:quotient-decay}: the torus has bounded diameter and the cigar
scalar curvature decays exponentially in its intrinsic distance. By
Proposition~\ref{prop:regular-lift-data} it is the quotient of a regular
circle lift, with the circle period chosen to make the transverse
K\"{a}hler class integral in the normalization of that proposition. The
universal cover of this quotient is the cigar times a flat complex line.
The circle bundle pulls back to a trivial bundle over this contractible
cover, and unwrapping its circle fibers gives the universal cover of the
Sasakian lift. Along the flat line over the cigar tip the scalar curvature
stays positive, so \eqref{eq:lifted-leaf-decay} fails by
\eqref{eq:lifted-leaf-distance}.

\section{Quasi-Regular Solitons and Orbifold Developability}

\label{sec:quasi-regular}

We work first on the quotient orbifold, without assuming developability. The
local soliton equation gives a parallel Ricci kernel. We then use this
splitting and the geometry of the potential to construct smooth covers in
transverse complex dimension two.

\subsection{Orbifold Quotients and the Ricci-Null Distribution}

\label{subsec:orbifold-geometry}

Let $\pi :M\rightarrow Q=M/\mathbb{S}^{1}$ be the complete effective K\"{a}%
hler orbifold quotient, and let $L$ be the generic Reeb period. The
descended potential satisfies, in every uniformizing chart, 
\begin{equation}
\func{Ric}_{h}+\func{Hess}_{h}f_{Q}=0.  \label{eq:orb-steady-equation}
\end{equation}%
In a chart centered at $q$, the circle orbibundle has the form 
\begin{equation*}
(\widehat{U}\times (\mathbb{R}/L\mathbb{Z}))/\Gamma _{q},\qquad \gamma
(z,t)=(\gamma z,t+\chi _{q}(\gamma )).
\end{equation*}%
An element of $\Gamma _{q}$ fixes a point above the chart center exactly
when its fiber character is trivial. Thus smoothness of $M$ forces $\chi
_{q} $ to be injective. Indeed, if $\gamma\in\ker\chi_q$, then $\gamma$
fixes every point of the circle fiber over the chart center, contradicting
freeness of the local action on the smooth total space unless $\gamma=1$.
Since every finite subgroup of the circle group is cyclic, the injective
character $\chi_q:\Gamma_q\to\mathbb{R}/L\mathbb{Z}$ also shows that $%
\Gamma_q$ is cyclic. This elementary observation is the only local
orbibundle fact needed below.

Now we first derive the constant rank theorem.

\begin{lemma}
\label{lem:orb-ricci-null} On a connected K\"{a}hler orbifold carrying a
gradient steady Ricci soliton with nonnegative holomorphic bisectional
curvature%
\begin{equation*}
B^{T}\geq 0,
\end{equation*}%
we have

\begin{enumerate}
\item The Ricci tensor has constant complex rank.

\item Its kernel is a parallel, $J$-invariant orbifold distribution.
\end{enumerate}
\end{lemma}

This is a local tensor assertion and does not assume the developability.

\begin{proof}
Set $f=f_{Q}$ and $S=\func{Ric}_{h}$. On tensors we use $L_{f}=\Delta
-\nabla _{\nabla f}$, where $\Delta $ is the rough real Laplacian. In a real
orthonormal frame, define 
\begin{equation*}
\mathcal{A}(S)(X,Y)=\sum_{a,b}S(e_{a},e_{b})\func{Rm}(X,e_{a},e_{b},Y).
\end{equation*}%
The local soliton flow gives 
\begin{equation*}
\partial _{t}S=\Delta S+2\mathcal{A}(S)-2S^{2},\qquad \mathcal{L}_{\nabla
f}S=\nabla _{\nabla f}S-2S^{2},
\end{equation*}%
and therefore 
\begin{equation}
L_{f}S=-2\mathcal{A}(S).  \label{eq:orb-static-ricci}
\end{equation}%
This is the covariant version of \cite[Lemma~2.1]{PetersenWylie}. Choose a
real unitary eigenframe $(e_{j},Je_{j})$ and write $%
r_{j}=S(e_{j},e_{j})=S(Je_{j},Je_{j})\geq 0$. The K\"{a}hler symmetries give 
\begin{equation}
\mathcal{A}(S)(X,X)=\sum_{j}r_{j}\func{Rm}(X,JX,Je_{j},e_{j})\geq 0.
\label{eq:orb-positive-reaction}
\end{equation}

(1) To prove constant rank, let $\sigma _{k}\geq 0$ be the sum of the $k$
smallest complex eigenvalues of $S$. Fix a point $p$ and extend the
corresponding eigenvectors to a local unitary frame normal at $p$. The
min--max characterization gives a smooth upper support $\tau
=\sum_{j=1}^{k}S(E_{j},E_{j})$ satisfying 
\begin{equation*}
\tau (p)=\sigma _{k}(p),\qquad (L_{f}\tau )(p)=-2\sum_{j=1}^{k}\mathcal{A}%
(S)(e_{j},e_{j})\leq 0.
\end{equation*}%
It follows that $L_{f}\sigma _{k}\leq 0$ in the barrier sense. The function $%
\sigma _{k}$ is invariant under every local isotropy group and hence is a
globally defined continuous orbifold function. Moreover, chart embeddings
are local isometries, so they identify the Ricci endomorphisms and carry the
above barrier supports to barrier supports. The strong minimum principle
therefore propagates across chart boundaries. Applying it in uniformizing
charts and using connectedness, each $\sigma _{k}$ is either identically
zero or everywhere positive. Thus the complex rank is constant.

(2) The kernel is therefore a smooth orbifold subbundle. For a local section 
$v$ of this kernel, differentiating $S(v,\cdot )=0$ gives 
\begin{equation*}
0=L_{f}(S(v,v))=-2\mathcal{A}(S)(v,v)-2\sum_{a}S(\nabla _{e_{a}}v,\nabla
_{e_{a}}v).
\end{equation*}%
Each term is nonpositive, so $\nabla _{e_{a}}v$ lies in the kernel for every 
$a$. Thus the kernel is parallel. Since $S$ is Hermitian, its kernel is $J$%
-invariant. The argument is equivariant in every chart, and the
distributions agree on overlaps because they are intrinsically defined.
\end{proof}

\subsection{Obstructions and Finite Smooth Covers}

\label{subsec:orbifold-covers}

In dimension five, we prove directly that the orbifold quotient of every
complete quasi-regular gradient steady Sasaki--Ricci soliton with
$B^{T}\geq 0$ is developable. We first establish results in arbitrary
dimension in the transversely flat case and in the positive-Ricci case
when the potential has a critical point.

\begin{proposition}
\label{prop:orb-flat-good} Suppose $M$ is complete and quasi-regular,
satisfies the steady equation and 
\begin{equation*}
B^{T}\geq 0.
\end{equation*}
If $R^{T}\equiv 0$, then 
\begin{equation*}
Q\cong \lbrack \mathbb{C}^{n}/\Lambda ]
\end{equation*}%
as a flat K\"{a}hler orbifold, where $\Lambda $ is a discrete group of
holomorphic Euclidean isometries acting properly with finite stabilizers. In
particular, $Q$ is developable. The same conclusion holds whenever $Q$ is
compact.
\end{proposition}

\begin{proof}
By Theorem~\ref{thm:conservation}, the universal cover of $M$ is the
Heisenberg model $H^{2n+1}\rightarrow \mathbb{C}^{n}$. Let $\Gamma $ be the
deck group and let $c_{s}$ denote central Reeb translation. If $L$ is the
generic period, then 
\begin{equation*}
\Gamma \cap \{c_{s}\}=\langle c_{L}\rangle .
\end{equation*}%
The intermediate quotient $P=H^{2n+1}/\langle c_{L}\rangle \rightarrow 
\mathbb{C}^{n}$ is a circle bundle. The action of $\Gamma /\langle
c_{L}\rangle $ projects faithfully to a group $\Lambda $ of holomorphic
Euclidean isometries. To check properness of this projected action, let $%
K\subset \mathbb{C}^{n}$ be compact. Its inverse image in $P$ is compact.
Proper discontinuity of the deck action on $P$ gives 
\begin{equation*}
\#\{\lambda \in \Lambda :\lambda K\cap K\neq \varnothing \}<\infty .
\end{equation*}%
Thus $Q=[\mathbb{C}^{n}/\Lambda ]$. If $Q$ is compact, orbifold integration
of $\Delta f_{Q}=-R_{h}$ first gives $R_{h}\equiv 0$, reducing to this case.
\end{proof}

\begin{proposition}
\label{prop:orb-compact-obstruction} Every K\"{a}hler orbifold satisfying %
\eqref{eq:orb-steady-equation} with a globally defined potential has 
\begin{equation}
c_{1}^{\mathrm{orb}}(TQ)_{\mathbb{R}}=0.  \label{eq:orb-c1-zero}
\end{equation}%
If in addition $B_{h}\geq 0$, then on every compact connected complex
suborbifold $S\subset Q$ of positive dimension, $f_{Q}|_{S}$ is constant and 
$\func{Ric}_{h}$ vanishes on $TS$. Thus no such $S$ exists when $\func{Ric}%
_{h}>0$. For a compact complex curve $S$, one also has $\chi _{\mathrm{orb}%
}(S)\leq 0$.
\end{proposition}

\begin{proof}
With $d^{c}u=-du\circ J$, the soliton equation gives 
\begin{equation*}
\rho _{h}=-\tfrac{1}{2}dd^{c}f_{Q},\qquad 2\pi c_{1}^{\mathrm{orb}}(TQ)_{%
\mathbb{R}}=[\rho _{h}]=0.
\end{equation*}%
Now put $F=-f_{Q}$. A complex suborbifold is minimal, so for compact $S$, 
\begin{equation*}
0=\int_{S}\Delta _{S}F\,dV=\int_{S}\func{tr}_{TS}\func{Ric}_{h}\,dV\geq 0.
\end{equation*}%
The nonnegative integrand must vanish. Hence $F|_{S}$ is harmonic and
therefore constant, while semidefiniteness gives $\func{Ric}_{h}|_{TS}=0$.
If $S$ is a complex curve, nonnegative bisectional curvature forces the
ambient holomorphic sectional curvature in its tangent direction to vanish.
The Gauss equation then gives $K_{S}\leq 0$, and orbifold Gauss--Bonnet
yields $2\pi \chi _{\mathrm{orb}}(S)=\int_{S}K_{S}\,dA\leq 0$.
\end{proof}

For $Q=\mathbb{C}\times\mathbb{CP}^1_{[m_1,m_2]}$, the compact slice has
trivial normal bundle because it is the product slice $\{0\}\times\mathbb{CP}%
^1_{[m_1,m_2]}$ and its normal orbibundle is the constant complex line $T_0%
\mathbb{C}$; hence 
\begin{equation*}
\int_{\{0\}\times\mathbb{CP}^1_{[m_1,m_2]}} c_1^{\mathrm{orb}%
}(TQ)=\frac1{m_1}+\frac1{m_2}>0.
\end{equation*}
Thus \eqref{eq:orb-c1-zero} excludes a global steady potential on these
orbifolds. In particular, any Sasakian construction whose transverse
quotient contains such a compact weighted-projective-line slice with
positive orbifold first Chern number cannot carry a global steady potential
satisfying \eqref{eq:orb-steady-equation}.

\begin{theorem}
\label{thm:orb-minimum-cover} Let $(Q^{2n},h,J,f_{Q})$ be a complete
noncompact K\"{a}hler orbifold satisfying \eqref{eq:orb-steady-equation}
with 
\begin{equation*}
\func{Ric}_{h}>0.
\end{equation*}
Suppose $f_{Q}$ has a critical point $o$. Then

\begin{enumerate}
\item there is a complete smooth simply connected K\"{a}hler gradient steady
soliton $(N,h_{N},J_{N},f_{N})$, diffeomorphic as a real manifold to $%
\mathbb{R}^{2n}$, and a finite group $\Gamma \subset U(n)$ such that 
\begin{equation}
(Q,h,J,f_{Q})\cong \lbrack (N,h_{N},J_{N},f_{N})/\Gamma ].
\label{eq:orb-finite-cover}
\end{equation}%
The group fixes the unique lifted critical point $\widetilde{o}$, is
isomorphic to the full local group at $o$, and acts by holomorphic
isometries preserving $f_{N}$.

\item for the quotient map $\Pi :N\rightarrow Q$, one has 
\begin{equation}
d_{h}(o,\Pi x)=d_{h_{N}}(\widetilde{o},x),\qquad R_{h_{N}}(x)=R_{h}(\Pi x),
\label{eq:orb-tip-distance}
\end{equation}%
Consequently, uniform scalar decay and either side of a linear decay
estimate pass to $N$. The required critical point exists whenever $%
R_{h}\rightarrow 0$ uniformly at infinity.
\end{enumerate}
\end{theorem}

\begin{proof}
The soliton identities hold in uniformizing charts and give $R_{h}+|\nabla
f_{Q}|^{2}=C_{0}$. If $R_{h}\rightarrow 0$ at infinity, it attains a
positive maximum. At that point, $0=\nabla R_{h}=2\func{Ric}_{h}(\nabla
f_{Q},\cdot )$ implies $\nabla f_{Q}=0$. It therefore suffices to construct
the cover from a given critical point.

Put $F=f_{Q}(o)-f_{Q}$. The geodesic argument in Proposition~\ref%
{prop:potential} applies to minimizing orbifold geodesics. Since $\func{Hess}%
F=\func{Ric}_{h}>0$, the point $o$ is the unique minimum, and for suitable $%
a,\mu >0$, 
\begin{equation}
\mu a(d_{h}(o,x)-a)\leq F(x)\leq \sqrt{C_{0}}\,d_{h}(o,x)\qquad
(d_{h}(o,x)\geq a).  \label{eq:orb-proper-potential}
\end{equation}%
Thus $F$ is proper and has range $[0,\infty )$.

Choose a uniformizing chart $(\widehat U,\Gamma)$ centered at $o$, and let $%
\widehat F$ be the lifted potential. Equivariant Morse coordinates identify
a sufficiently small lifted sublevel $\widehat\Omega_{2\varepsilon}=\{%
\widehat F<2\varepsilon\}$ with a ball and $\widehat\Sigma_\varepsilon=\{%
\widehat F=\varepsilon\}$ with a sphere. The isotropy representation at the
center is faithful and unitary. Properness and uniqueness of the minimum
allow us to choose $\varepsilon$ so that the entire sublevel $%
\{F<2\varepsilon\}$ is contained in this chart.

We extend this local covering along the level sets. On $Q\setminus\{o\}$,
the orbifold vector field $V=\nabla F/|\nabla F|^2$ satisfies $VF=1$. Its
flow is defined between any two positive levels, since the intervening slab
is compact and contains no critical point. It therefore identifies 
\begin{equation*}
Q\setminus\{o\}\cong\Sigma_\varepsilon\times(0,\infty).
\end{equation*}
Glue 
\begin{equation*}
N=\widehat\Omega_{2\varepsilon}\ \cup\ \bigl(\widehat\Sigma_\varepsilon%
\times(\varepsilon,\infty)\bigr)
\end{equation*}
along $\widehat\Sigma_\varepsilon\times(\varepsilon,2\varepsilon)$ using the
lifted $V$-flow. The local flows are equivariant and conjugate the
stabilizer actions. Thus the chart projections agree on the overlap and
define an orbifold covering $\Pi:N\to Q$ from a smooth manifold. The group $%
\Gamma$ acts on both pieces and on their gluing, and $Q=[N/\Gamma]$. The
construction attaches an infinite collar to a ball, giving $N\cong\mathbb{R}%
^{2n}$; it is equivariant with respect to the linear isotropy action.

The gluing is performed on the smooth lifted ball and sphere collar, so
$N$ is a smooth manifold. For $x\in N$, let $\Gamma_x$ denote its
stabilizer and choose a sufficiently small $\Gamma_x$-invariant
neighborhood $U_x$, disjoint from its translates by elements outside
$\Gamma_x$. The quotient map is locally $U_x/\{1\}\to U_x/\Gamma_x$,
corresponding to the subgroup inclusion $\{1\}\subset\Gamma_x$; see
\cite[Section~2.1, p.~107]{KleinerLott2014}. Thus the covering is understood
in the orbifold sense; no freeness of $\Gamma$ away from $\widetilde o$
is required.

The chart pullbacks of $h,J,f_Q$ now define smooth tensors $h_N,J_N,f_N$
satisfying the K\"ahler and steady equations. To prove completeness, let $%
(x_j)$ be an $h_N$-Cauchy sequence. Its projection converges in $Q$. The
finite covering has finitely many precompact lifted chart neighborhoods over
a small neighborhood of that limit, so $(x_j)$ has a convergent subsequence.
The Cauchy property then gives convergence of the entire sequence.

Finally, $\Pi^{-1}(o)=\{\widetilde o\}$ and every element of $\Gamma$ fixes $%
\widetilde o$. The distance formula for a finite isometric quotient gives 
\begin{equation*}
d_h(o,\Pi x)=\min_{\gamma\in\Gamma} d_{h_N}(\widetilde o,\gamma
x)=d_{h_N}(\widetilde o,x).
\end{equation*}
Scalar curvature pulls back locally, proving \eqref{eq:orb-tip-distance} and
the assertions about decay.
\end{proof}

\begin{corollary}
\label{cor:orb-decay-cyclic} Let $M^{2n+1}$, $n\geq 1$, be a complete
quasi-regular gradient steady Sasaki-Ricci soliton with $B^{T}\geq 0$.
Suppose that its orbifold quotient $Q^{2n}$ has $\func{Ric}_{h}>0$ and $%
R_{h}\rightarrow 0$ uniformly at infinity. Then

\begin{enumerate}
\item $Q$ is developable and is a finite cyclic quotient of a complete
simply connected irreducible K\"{a}hler steady manifold $N^{2n}$ with
positive Ricci curvature, nonnegative bisectional curvature and uniform
scalar decay.

\item Moreover, $N\cong \mathbb{R}^{2n}$ as a smooth real manifold and, for
the generic period $L$, the total space has a presentation 
\begin{equation}
M\cong \bigl(N\times (\mathbb{R}/L\mathbb{Z})\bigr)/\Gamma ,\qquad \gamma
\cdot (x,t)=(\gamma x,t+\chi (\gamma )),  \label{eq:orb-total-cyclic}
\end{equation}%
where $\Gamma $ is finite cyclic, fixes the tip of $N$, and $\chi :\Gamma
\rightarrow \mathbb{R}/L\mathbb{Z}$ is faithful. The lifted contact form is $%
dt+A$, where $dA=2\omega _{N}$ and $A$ may be chosen $\Gamma $-invariant.
The Reeb orbit over $x$ has period $L/|\Gamma _{x}|$.
\end{enumerate}
\end{corollary}

\begin{proof}
By Theorem~\ref{thm:orb-minimum-cover}, $Q=[N/\Gamma ]$ and scalar decay
lifts to $N$. The group $\Gamma $ is the full isotropy group at the critical
point. The faithful fiber character therefore makes it cyclic. The splitting
theorem and Proposition~\ref{prop:product-obstruction} imply that $N$ is
irreducible.

The pulled-back circle bundle over $N\cong \mathbb{R}^{2n}$ is trivial.
Equivariant Morse coordinates near the minimum and the level flow outside
give an equivariant contraction to the critical point. Parallel transport in
the connection along these radial paths gives an equivariant trivialization,
in which 
\begin{equation*}
\gamma (x,t)=(\gamma x,t+\chi (\gamma )),\qquad \eta =dt+A,\qquad \gamma
^{\ast }A=A,\qquad dA=2\omega _{N}.
\end{equation*}%
Freeness of the total-space action over the critical point makes $\chi $
faithful. Over $x$, the stabilizer acts on the circle by a subgroup of order 
$|\Gamma _{x}|$, so the Reeb period is $L/|\Gamma _{x}|$.
\end{proof}

Conversely, the data in \eqref{eq:orb-total-cyclic}, with faithful $\chi$,
give a smooth complete lift. An invariant primitive is obtained by averaging
over $\Gamma$.

\subsection{Developability in Dimension Five}

\label{subsec:orbifold-five}

\begin{proposition}
\label{prop:orb-positive-surface-good} Let $(Q^{4},h,J,f_{Q})$ be a complete
K\"{a}hler orbifold satisfying \eqref{eq:orb-steady-equation}, with 
\begin{equation*}
B_{h}\geq 0\text{ \ \ and \ \ }\func{Ric}_{h}>0.
\end{equation*}
Then

\begin{enumerate}
\item $Q$ is developable. If $Q$ has any orbifold singularity, its potential
has a critical point, and $Q$ has the finite-quotient description of Theorem~%
\ref{thm:orb-minimum-cover}.

\item For the quotient of a smooth quasi-regular Sasakian five-manifold, the
finite group in this description is cyclic.
\end{enumerate}
\end{proposition}

\begin{proof}
There is nothing to prove if $Q$ is smooth. Otherwise put $F=-f_{Q}$ and
suppose $\nabla F$ never vanishes. At a singular point its nonzero gradient
is fixed by $\Gamma _{q}\subset U(2)$. Effectiveness forces the common fixed
space to be a complex line, while the normal action is a faithful finite
subgroup of $U(1)$. Thus the singular locus is locally a totally geodesic
complex disk with constant cyclic isotropy.

Choose a connected component $S$ of the singular locus. We record why the
local fixed disks glue globally. On an overlap of two effective uniformizing
charts, a chart embedding conjugates the stabilizer of a point into the
stabilizer in the larger chart and carries its fixed-point set to the
corresponding fixed-point set. The no-critical-point assumption is used at
every singular point : $\nabla F(q)\neq 0$ is fixed by the full local
isotropy group, so its complex span is contained in the fixed space.
Effectiveness excludes a two-dimensional fixed space for a nontrivial
stabilizer. Hence the fixed space has complex dimension exactly one at every
point of the singular locus. A chart embedding is equivariant for the
injective homomorphism of stabilizers, and its differential maps this
one-dimensional fixed space onto the one-dimensional fixed space in the
target chart. Therefore the fixed disks agree on overlaps. Hence $S$ is a
smooth embedded complex curve. There is no branching, and closedness of the
orbifold singular set rules out endpoints. The fixed disks are totally
geodesic, so $S$ is totally geodesic. A Cauchy sequence in $S$ is Cauchy in
the complete orbifold $Q$ and converges to the closed set $S$. Thus the
induced metric on $S$ is complete.

The potential and its gradient are invariant under every local isotropy
group. At a point of $S$ the invariant tangent space is exactly $TS$, so $%
\nabla F$ is tangent to $S$ globally. Restriction of the soliton equation
gives 
\begin{equation*}
\func{Hess}_{S}F=\psi h_{S},\qquad \psi >0,\qquad K_{S}\geq 0.
\end{equation*}%
Let $w=|\nabla _{S}F|>0$ and $\nu =\nabla _{S}F/w$. From $\nabla
_{S}^{2}F=\psi h_{S}$ one has $\nabla _{\nu }\nabla F=\psi \nu $, hence $\nu
(w)=\psi $ and $\nabla _{\nu }\nu =0$. Because $S$ is complete, the
unit-speed geodesics tangent to $\nu $ exist for all time. Along such a
geodesic $\gamma $ choose a parallel unit field $e\perp \nu $. The Ricci
commutation identity applied to $F$ gives 
\begin{equation*}
(\nabla _{\nu }\nabla ^{2}F)(e,e)-(\nabla _{e}\nabla ^{2}F)(\nu ,e)=\langle 
\mathcal{R}(\nu ,e)\nabla F,e\rangle =-K_{S}w.
\end{equation*}%
Since $\nabla ^{2}F=\psi h_{S}$, the left-hand side is $\nu (\psi )$.
Consequently 
\begin{equation*}
w^{\prime }=\psi >0,\qquad \psi ^{\prime }=-K_{S}w,\qquad w^{\prime \prime
}=-K_{S}w\leq 0.
\end{equation*}%
Concavity implies $w(s)\leq w(0)+s\psi (0)$ for $s<0$, contradicting $w>0$.
Thus $F$ has a critical point, and Theorem~\ref{thm:orb-minimum-cover}
applies. In the Sasakian setting, the faithful fiber character makes the
finite group cyclic, exactly as in Corollary~\ref{cor:orb-decay-cyclic}.
Note that a finite group admitting an injective character into $S^{1}$ is
cyclic.
\end{proof}

\begin{proposition}
\label{prop:orb-rank-one-good} Let $(Q^{4},h,J,f_{Q})$ be a complete K\"{a}%
hler orbifold satisfying \eqref{eq:orb-steady-equation}, with $B_{h}\geq 0$
and complex Ricci rank one. Then its universal orbifold cover is the smooth K%
\"{a}hler product 
\begin{equation*}
(\mathbb{C},ag_{\mathrm{cig}})\times (\mathbb{C},g_{\mathrm{flat}}),\qquad
a>0.
\end{equation*}%
In particular, $Q$ is developable.
\end{proposition}

\begin{proof}
On the complete universal orbifold cover $\widehat Q$, the distributions $%
\mathcal{E}=\ker\func{Ric}$ and $\mathcal{P}=\mathcal{E}^\perp$ are parallel
and $J$-invariant by Lemma~\ref{lem:orb-ricci-null}. At a regular point they
give an orthogonal tangent splitting invariant under holonomy. The orbifold
de Rham theorem \cite[Lemma~2.19]{KleinerLott2014}, which applies in every
dimension, therefore gives an isometric splitting 
\begin{equation*}
\widehat Q=E\times\Sigma.
\end{equation*}
The factors are complete because their slices are closed and totally
geodesic in the complete product. They are simply connected in the orbifold
sense, since the orbifold fundamental group of a product is the product of
those of its factors and $\pi_1^{\mathrm{orb}}(\widehat Q)=0$. The parallel
complex structure preserves the two factors, so this is a K\"ahler product
of orbifold curves. The factor $E$ is flat, and $\Sigma$ has positive
Gaussian curvature. The mixed Hessian equation splits the lifted potential
into a sum of functions on the two factors; each factor therefore satisfies
the steady soliton equation.

The flat connection on $E$ has holonomy factoring through $\pi_1^{\mathrm{orb%
}}(E)=0$. In particular, all local isotropy representations are trivial.
Effectiveness therefore makes $E$ a smooth complete simply connected flat
curve, hence $E=\mathbb{C}$.

Suppose that $\Sigma$ has a singular point. Its nontrivial local group is a
subgroup of $U(1)$ with no nonzero fixed vector, so the lifted potential
gradient vanishes at that point. Theorem~\ref{thm:orb-minimum-cover} would
express $\Sigma$ as a finite quotient of a smooth simply connected steady
surface. Because the smooth cover is simply connected, its deck group is $%
\pi_1^{\mathrm{orb}}(\Sigma)$, which is trivial. This contradicts the
singularity. Hence $\Sigma$ is smooth, and Lemma~\ref{lem:surface}
identifies it with a scaled cigar.
\end{proof}

\begin{theorem}
\label{thm:quasi-all-good} Let $(M^{5},g,\eta ,\xi ,\Phi ,f)$ be a complete
quasi-regular gradient steady Sasaki-Ricci soliton with 
\begin{equation*}
B^{T}\geq 0.
\end{equation*}%
Its effective K\"{a}hler orbifold quotient is developable, and its smooth
universal orbifold cover has one of the four forms in Theorem~\ref{thm:five}.
\end{theorem}

\begin{proof}
Lemma~\ref{lem:orb-ricci-null} gives constant complex Ricci rank $r\in
\{0,1,2\}$. Propositions~\ref{prop:orb-flat-good}, \ref%
{prop:orb-rank-one-good} and~\ref{prop:orb-positive-surface-good} prove
developability for $r=0,1,2$, respectively. Apply Theorem~\ref{thm:five} to
the resulting smooth complete cover.
\end{proof}

\begin{remark}
The total space is recovered by pulling the circle orbibundle back to $N$.
This gives a principal circle bundle $\widetilde{M}\rightarrow N$ and an
ordinary covering $\widetilde{M}\rightarrow M$ with deck group $\Gamma =\pi
_{1}^{\mathrm{orb}}(Q)$. In terms of data on $Q$, the lift is specified by
an integral class 
\begin{equation}
e\in H_{\mathrm{orb}}^{2}(Q;\mathbb{Z}),\qquad e_{\mathbb{R}}=\frac{2}{L}%
[\omega _{h}],  \label{eq:steady-orb-euler}
\end{equation}%
with faithful local fiber characters, together with a connection of
curvature $2\omega _{h}$. No finiteness of $\Gamma$ is asserted at this
level of generality. The finite cyclic conclusion of Corollary~\ref%
{cor:orb-decay-cyclic} uses the critical point of a potential with positive
definite Hessian.
\end{remark}

\section{Irregular Solitons and Complete Transverse Developments}

\label{sec:irregular}

We now construct transverse developments directly from the foliation. The
argument has two parts: a local curvature identity produces a critical point
when the transverse Ricci tensor is positive, and a neighborhood of its
critical leaf extends to a global model by the potential flow. Complex Ricci
rank one in dimension five is treated separately by continuation of the
local cigar-line geometry. Under their common hypotheses, the
monotonicity, critical-point, and positive-Ricci development arguments
are shared with \cite[Lemma~3.1 and Theorems~3.1--3.2]{Rigidity2026};
we give their proofs here for the
dimension-five construction and the subsequent suspension classification.

\subsection{Ricci Monotonicity and Critical Points}

\label{subsec:critical-points}

Put 
\begin{equation*}
F=-f\text{ \ \ \textrm{and} \ \ }V:=\nabla ^{T}F=\nabla ^{g}F
\end{equation*}%
and 
\begin{equation*}
S=\func{Ric}^{T}=\func{Hess}^{T}F
\end{equation*}

\begin{lemma}
\label{lem:transverse-ricci-monotonicity} For a gradient steady Sasaki-Ricci
soliton, write $W=V^{1,0}=\tfrac{1}{2}(V-iJ^{T}V)$. In transverse K\"{a}hler
coordinates,

\begin{enumerate}
\item 
\begin{equation}
(\nabla _{V}^{T}S)_{i\bar{j}}=-2R_{i\bar{j}k\bar{\ell}}^{T}W^{k}\overline{%
W^{\ell }}.  \label{eq:transverse-ricci-monotonicity}
\end{equation}

\item If 
\begin{equation*}
B^{T}\geq 0,
\end{equation*}%
then $\nabla _{V}^{T}S$ is a nonpositive Hermitian form, equivalently a
nonpositive real symmetric form.
\end{enumerate}
\end{lemma}

\begin{proof}
Work in a transverse K\"{a}hler chart and suppress the superscript $T$.
Since $\func{Hess}F=S$ is Hermitian, 
\begin{equation}
F_{ik}=0,\qquad S_{i\bar{j}}=F_{i\bar{j}},\qquad W^{k}=h^{k\bar{\ell}}F_{%
\bar{\ell}}.  \label{eq:critical-hessian-components}
\end{equation}%
Choose holomorphic normal coordinates at a point and differentiate $%
F_{ik}=\partial _{i}\partial _{k}F-\Gamma _{ik}^{\ell }F_{\ell }=0$ : 
\begin{equation}
0=\partial _{\bar{j}}F_{ik}=\nabla _{k}S_{i\bar{j}}+R_{i\bar{j}k}{}^{\ell
}F_{\ell }.  \label{eq:critical-commutation}
\end{equation}%
Thus 
\begin{equation}
\nabla _{k}S_{i\bar{j}}=-R_{i\bar{j}k\bar{\ell}}\overline{W^{\ell }},\qquad
\nabla _{\bar{k}}S_{i\bar{j}}=-R_{i\bar{j}\ell \bar{k}}W^{\ell }.
\label{eq:critical-ricci-derivative}
\end{equation}%
Contracting with $W^{k}$ and $\overline{W^{k}}$ proves %
\eqref{eq:transverse-ricci-monotonicity}, since $V=W+\overline{W}$. For
every $(1,0)$ vector $U$, 
\begin{equation*}
(\nabla _{V}S)(U,\overline{U})=-2\func{Rm}(U,\overline{U},W,\overline{W}%
)\leq 0.
\end{equation*}%
The identities are intrinsic and hence agree across transverse charts.
\end{proof}

Identity~\eqref{eq:critical-ricci-derivative} was derived above directly
from the soliton equation and curvature commutation; no result from \cite%
{ChangHanLinWu2026} is used in its proof. We cite that preprint only for
comparison, since the same differentiated identity appears there in
equation~(3.3).

\begin{theorem}
\label{thm:positive-ricci-critical} Let $(M^{2n+1},g,\eta ,\xi ,\Phi ,f)$ be
a complete gradient steady Sasaki-Ricci soliton with 
\begin{equation*}
B^{T}\geq 0\text{ \ \ \textrm{and }\ \ }\func{Ric}^{T}>0.
\end{equation*}%
Then $f$ has a transverse critical point. At every such point $q$, 
\begin{equation}
R^{T}(q)=C_{0}=\max_{M}R^{T}.  \label{eq:critical-scalar-maximum}
\end{equation}%
The conclusion holds for every Reeb type without assuming transverse
developability.
\end{theorem}

\begin{proof}
By conservation, $R^{T}+|V|^{2}=C_{0}$ and $|V|\leq \sqrt{C_{0}}$. The
bounded field $-V$ is complete on the complete total manifold. Fix $p\in M$.
If $V_{p}=0$, the assertion follows immediately from conservation.
Otherwise, let 
\begin{equation}
\gamma ^{\prime }(t)=-V_{\gamma (t)},\qquad \gamma (0)=p,\qquad t\geq 0.
\label{eq:critical-gradient-trajectory}
\end{equation}%
Let $c>0$ be the least eigenvalue of $S_{p}$ with respect to $g_{p}^{T}$.
Parallel transport in the transverse connection is defined along all of $%
\gamma $. For a parallel horizontal field $Y$, the preceding lemma gives 
\begin{equation*}
\frac{d}{dt}S(Y,Y)=(\nabla _{-V}^{T}S)(Y,Y)\geq 0,\qquad |Y(t)|=|Y(0)|.
\end{equation*}%
Consequently 
\begin{equation}
S_{\gamma (t)}\geq c\,g_{\gamma (t)}^{T}.
\label{eq:critical-orbit-lower-bound}
\end{equation}%
For $u(t)=|V|^{2}(\gamma (t))$, 
\begin{equation}
u^{\prime }=-2S(V,V)\leq -2cu,  \label{eq:critical-speed-ode}
\end{equation}%
and hence 
\begin{equation}
|\gamma ^{\prime }(t)|\leq |V_{p}|e^{-ct},\qquad \func{Length}(\gamma
|_{[t,\infty )})\leq \frac{|V_{p}|}{c}e^{-ct}.
\label{eq:critical-finite-length}
\end{equation}%
The finite tail length makes $\gamma (t)$ Cauchy as $t\rightarrow \infty $.
Completeness gives a limit $q\in M$, and the speed estimate implies $V_{q}=0$%
. Conservation then gives $R^{T}(q)=C_{0}=\max_{M}R^{T}$.
\end{proof}

For a complete smooth K\"ahler soliton, critical-point existence is known
from \cite[Theorem~1.2]{DengZhu2015}. The proof above uses transverse
parallel transport along a curve in the total manifold and therefore applies
before a global transverse model is available. We next use the resulting
critical leaf to construct that model.

\subsection{Construction of Complete Transverse Developments}

\label{subsec:development-constructions}

\begin{proposition}
\label{prop:rank-one-development} Let $M^{5}$ be a complete gradient steady
Sasaki-Ricci soliton with 
\begin{equation*}
B^{T}\geq 0.
\end{equation*}%
If $\func{Ric}^{T}$ has complex rank one at some point, then

\begin{enumerate}
\item its rank is one everywhere.

\item $M$ admits a complete smooth transverse development 
\begin{equation*}
X=(\mathbb{C},ag_{\mathrm{cig}})\times (\mathbb{C},g_{\mathrm{flat}}),\qquad
a>0.
\end{equation*}
\end{enumerate}
\end{proposition}

\begin{proof}
We first identify a common local model. The proof of Lemma~\ref%
{lem:orb-ricci-null} applies in transverse charts and gives parallel complex
lines $E=\ker \func{Ric}^{T}$ and $P=E^{\perp }$. The projection $V_{0}=%
\func{pr}_{E}\nabla ^{T}f$ is transversely parallel. Hence $b=|V_{0}|^{2}$
is constant. Local product coordinates split the metric into a flat line and
a surface of positive Gaussian curvature, and split the potential into an
affine flat summand and a surface summand $f_{+}$.

Put $F_{+}=-f_{+}$. On the surface, 
\begin{equation*}
\func{Hess}F_{+}=Kh_{+},\qquad dK=-K\,dF_{+},\qquad 2K+|\nabla
F_{+}|^{2}=C_{+}:=C_{0}-b>0.
\end{equation*}%
Away from a critical point, use a flow coordinate $\tau $ for the Killing
field $J\nabla F_{+}$. Since $Ke^{F_{+}}=\alpha $ is a positive constant,
the metric is 
\begin{equation*}
h_{+}=\frac{dF_{+}^{2}}{C_{+}-2\alpha e^{-F_{+}}}+(C_{+}-2\alpha
e^{-F_{+}})d\tau ^{2}.
\end{equation*}%
The substitution 
\begin{equation*}
r^{2}=\frac{C_{+}}{2K}-1,\qquad \theta =\frac{C_{+}}{2}\tau
\end{equation*}%
gives 
\begin{equation*}
h_{+}=(4/C_{+})(dr^{2}+r^{2}d\theta ^{2})/(1+r^{2}).
\end{equation*}
At a critical point, the local smooth polar calculation in the proof of
Lemma~\ref{lem:surface} gives the cigar tip; it does not require
completeness of the transverse chart. Thus every transverse chart maps into
the fixed complete model 
\begin{equation*}
X=(\mathbb{C},(4/C_{+})g_{\mathrm{cig}})\times \mathbb{C}.
\end{equation*}

We next justify the continuation step used to globalize this local model.
Every holomorphic local isometry $\phi:U\to V$ between connected open
subsets of $X$ preserves the Ricci tensor. Hence $d\phi(\ker\func{Ric})=\ker%
\func{Ric}$ and likewise for its orthogonal complement. These are precisely
the tangent distributions of the flat and cigar factors, respectively. Thus,
in product coordinates, $\phi(z,w)=(\phi_1(z),\phi_2(w))$ locally; no mixing
of the two complex curves is possible.

The cigar scalar curvature has a unique maximum at its tip and is a strictly
decreasing function of the cigar distance from that tip. A local cigar
isometry preserves scalar curvature, hence preserves that distance function
on its domain. In polar coordinates it therefore preserves the radial
coordinate; holomorphicity then forces the angular change to be constant.
Thus $\phi_1$ is the restriction of a unique global rotation of the cigar.
The flat component $\phi_2$ is a holomorphic Euclidean local isometry, so $%
\phi_2(w)=e^{i\theta}w+b$ and also extends uniquely. Consequently every
holomorphic local isometry of a connected open subset of $X$ extends
uniquely to an element of $\func{Isom}_{\mathrm{hol}}(X)$. Uniqueness
implies that the extended maps satisfy exactly the same cocycle identities
as the original transitions.

Starting from one transverse chart, continue its model map along a path in $%
M $ by these global extensions. Homotopic paths with fixed endpoints give
the same continuation: subdivide a homotopy into chart rectangles and use
the cocycle identity on each rectangle. Hence continuation is single-valued
on the universal cover $q:\widehat{M}\to M$ and defines a transverse local
submersion $D:\widehat{M}\to X$.

We verify completeness and connectedness of the development. A horizontal
lift of a piecewise smooth path has the same length as the base path.
Completeness of $\widehat{M}$ prevents finite-time escape, so every base
path on a compact interval lifts from any chosen initial point. In
particular, $D$ is surjective. Transport along radial paths in a strongly
convex ball gives local trivializations. The homotopy sequence, together
with contractibility of $X$ and simple connectivity of $\widehat{M}$, shows
that each fiber is connected and simply connected. Being one-dimensional and
tangent to the lifted Reeb field, it is a Reeb line.

Consequently $q^{\ast }f$ is constant on each fiber and descends to a smooth
function $f_{X}$. The transverse metric and complex structure descend with
it, and 
\begin{equation*}
\func{Ric}_{X}+\func{Hess}_{X}f_{X}=0.
\end{equation*}
The potential is the sum of the cigar potential and an affine function on
the flat factor. Thus all conditions in Definition~\ref{def:development}
hold.
\end{proof}

\begin{proposition}
\label{prop:critical-leaf-development} Let $(M^{2n+1},g,\eta ,\xi ,\Phi ,f)$%
, $n\geq 1$, be a complete gradient steady Sasaki-Ricci soliton with 
\begin{equation*}
\func{Ric}^{T}>0.
\end{equation*}
Suppose that the closure $H$ of the Reeb flow in $\func{Isom}(M,g)$ is
compact and that $\nabla ^{T}f(p)=0$ at some point $p$. Then

\begin{enumerate}
\item the critical set of $f$ is a single compact Reeb orbit $L$. The
function $F=f(p)-f$ is nonnegative and proper, and there are constants $c>0$
and $C<\infty $ such that 
\begin{equation*}
F(x)\geq c\,d_{g}(L,x)-C.
\end{equation*}

\item there is also a complete smooth transverse development $%
(N^{2n},h,J,f_{N})$ as in Definition~\ref{def:development}. The manifold $N$
is biholomorphic to $\mathbb{C}^{n}$, has positive Ricci curvature, and has
a unique critical point of $f_{N}$. The development can be taken on the
Riemannian universal cover of $M$. Its fibers are the lifted Reeb lines.
\end{enumerate}
\end{proposition}

\begin{proof}
We first identify the critical leaf and prove properness of the potential.
The Reeb closure preserves $f$, so $Hp\subset \func{Crit}(f)$. Since $\func{%
Hess}^{T}F=\func{Ric}^{T}>0$, the critical set in a transverse chart is an
isolated point and upstairs it is a local Reeb plaque. Hence $Hp$ is
one-dimensional. Compactness and the nonvanishing Reeb field make it a
circle $L$, with minimal period $\ell >0$.

For $x\notin L$, choose a minimizing unit-speed geodesic $\gamma
:[0,r]\rightarrow M$ from the compact submanifold $L$ to $x$, where $%
r=d_g(L,x)$. The first-variation formula for distance to a submanifold gives 
$\dot\gamma(0)\perp T_{\gamma(0)}L$, and since $T L=\mathbb{R}\xi$ this
means $\eta(\dot\gamma(0))=0$. Along any geodesic, 
\begin{equation*}
\frac d{ds}g(\dot\gamma,\xi)=g(\dot\gamma,\nabla_{\dot\gamma}\xi)=0,
\end{equation*}
because $\xi$ is Killing and therefore $g(Y,\nabla_Y\xi)=0$ for every $Y$.
Hence $\gamma$ remains horizontal. Thus 
\begin{equation*}
(F\circ \gamma )^{\prime }(0)=0,\qquad (F\circ \gamma )^{\prime \prime }=%
\func{Ric}^{T}(\dot{\gamma},\dot{\gamma})>0.
\end{equation*}%
In particular, $F(x)>0$ and $dF_{x}(\dot{\gamma}(r))>0$, so $L$ is the
entire critical set. Choose $\delta ,\mu >0$ with $\func{Ric}^{T}\geq \mu
g^{T}$ on the compact closed $\delta $-neighborhood of $L$. Integrating
along the same geodesics gives 
\begin{equation*}
F(x)\geq \mu \delta \bigl(d_{g}(L,x)-\delta \bigr)\qquad (d_{g}(L,x)\geq
\delta ).
\end{equation*}%
This proves the growth estimate and properness.

We next describe a neighborhood of $L$. For a sufficiently small normal tube 
$U$, the map 
\begin{equation*}
q_{0}:B_{\varepsilon }(\mathcal{D}_{p})\times \mathbb{R}\rightarrow U,\qquad
q_{0}(v,t)=\varphi _{t}(\exp _{p}v)
\end{equation*}%
is a covering. It is the normal disk bundle pulled back along $\mathbb{R}%
\rightarrow L$ and trivialized by the Reeb isometries. Its deck generator is 
\begin{equation*}
\sigma (v,t)=(Av,t+\ell ),\qquad A=d\varphi _{-\ell }|_{\mathcal{D}_{p}}\in
U(n).
\end{equation*}%
All pulled-back tensors are independent of $t$ and define a transverse K\"{a}%
hler metric on the disk; the exponential section itself need not be
horizontal.

Choose $a>0$ so that $\{F<2a\}\subset U$, and write $q_{0}^{\ast }F=F_{0}(v)$%
. The function $F_0$ is invariant under the compact closure of $\langle
A\rangle\subset U(n)$ and has positive-definite Hessian at the origin. The
equivariant Morse lemma therefore identifies $\Omega _{2a}=\{F_{0}<2a\}$
with a ball and $S_{a}=\{F_{0}=a\}$ with a sphere, compatibly with the
action of $A$. In particular, both are $A$-invariant, and 
\begin{equation*}
\{F<2a\}\cong (\Omega _{2a}\times \mathbb{R})/\langle \sigma \rangle ,\qquad
\{F=a\}\cong (S_{a}\times \mathbb{R})/\langle \sigma \rangle .
\end{equation*}

To extend this model, put $V=\nabla F/|\nabla F|^{2}$ on $M\setminus L$.
Then 
\begin{equation*}
VF=1,\qquad \lbrack V,\xi ]=0.
\end{equation*}%
The flow exists between any two positive levels because the intervening slab
is compact. Moreover, $F$ is unbounded: otherwise properness would make $M$
compact, and a maximum away from $L$ would be another critical point. Hence $%
F(M)=[0,\infty )$, and the level flow identifies $\{F>a\}$ with $%
\{F=a\}\times (a,\infty )$.

Glue the covers $\Omega _{2a}\times \mathbb{R}$ and $S_{a}\times \mathbb{R}%
\times (a,\infty )$ over $a<F<2a$ by the lifted flow. The relation $[V,\xi
]=0$ makes the transition equivariant under fiber translations. It also
intertwines the lifted deck action; the fiber coordinate may change by a
smooth function of the base point. Thus the gluing defines a covering $q:%
\widehat{M}\rightarrow M$ and a principal $\mathbb{R}$-bundle 
\begin{equation*}
D:\widehat{M}\longrightarrow N=\Omega _{2a}\ \cup \ \bigl(S_{a}\times
(a,\infty )\bigr)\cong \mathbb{R}^{2n}.
\end{equation*}%
Since $N$ is contractible, the bundle is trivial and $\widehat{M}$ is simply
connected. Therefore $q$ is the Riemannian universal covering. The
transverse tensors descend to $h,J,f_{N}$ on $N$, with $\func{Ric}_{h}+\func{%
Hess}_{h}f_{N}=0$ and $\func{Ric}_{h}>0$.

It remains to check completeness of $h$. On the two covering pieces, the
transverse tensors are pulled back from $M$. Since the overlap maps
intertwine the covering maps to $M$, these pulled-back tensors agree on
overlaps. Thus the descended metric satisfies $q^*g^T=D^*h$, and $D$ is
globally a Riemannian submersion with horizontal distribution
$\ker(q^*\eta)$. The covering metric $q^*g$ is complete
because $q$ is a Riemannian covering of the complete manifold $M$. Given $%
(y,v)\in TN$, choose $\hat x\in D^{-1}(y)$ and the unique horizontal $\hat v$
with $dD(\hat v)=v$. The $q^*g$-geodesic with initial data $(\hat x,\hat v)$
exists for all time. For a Riemannian submersion whose fibers are the
geodesic Reeb lines, a geodesic initially orthogonal to a fiber remains
horizontal; equivalently this follows from conservation of $%
g(\dot\gamma,\xi) $ since $\xi$ is Killing. Its projection is therefore the $%
h$-geodesic with initial data $(y,v)$, defined for all time. Hence $h$ is
geodesically complete. The fibers of $D$ are the lifted Reeb lines, and the
unique critical leaf projects to the unique critical point of $f_{N}$.
Conservation gives $R_h+|\nabla f_N|^2=C_0$, so $R_h$ attains its maximum at
that point. Since $N$ is complete and noncompact with positive Ricci
curvature, \cite[Theorem~1]{ChauTam2005} shows that $N$ is biholomorphic to $%
\mathbb{C}^n$. This proves the proposition.
\end{proof}

\begin{remark}
If $\func{Ric}^{T}>0$, any maximum point of $R^{T}$ is transversely
critical, since $\nabla ^{T}R^{T}=2\func{Ric}^{T}(\nabla ^{T}f,\cdot )$.
Thus the proposition also applies when $H$ is compact and $R^{T}\rightarrow
0 $ with distance from a compact $H$-orbit, without assuming $B^{T}\geq 0$.
Under $B^{T}\geq 0$, Theorem~\ref{thm:positive-ricci-critical} supplies the
critical point without decay. A proper Reeb line action already has a
complete smooth quotient.
\end{remark}

\begin{theorem}
\label{thm:five-all-development} Every complete five-dimensional gradient
steady Sasaki-Ricci soliton with%
\begin{equation*}
B^{T}\geq 0
\end{equation*}
admits a complete smooth transverse development in the sense of Definition~%
\ref{def:development}.
\end{theorem}

\begin{proof}
Regular Reeb circle actions and proper Reeb line actions have complete
smooth quotients. The quasi-regular case follows from Theorem~\ref%
{thm:quasi-all-good}. For an irregular flow, $H$ is a compact torus and the
complex Ricci rank is constant. In rank zero, Theorem~\ref{thm:conservation}
gives the Heisenberg development. Rank one is Proposition~\ref%
{prop:rank-one-development}. In rank two, Theorem~\ref%
{thm:positive-ricci-critical} gives a critical point, and Proposition~\ref%
{prop:critical-leaf-development} supplies the development.
\end{proof}

Theorem~\ref{thm:intro-five} now follows from this development theorem and
the classification in Theorem~\ref{thm:five}.

\subsection{Suspensions and Steady Approximations}

\label{subsec:suspension-classification}

In the positive Ricci case, we determine the global suspension data and
construct quasi-regular steady approximations.

Let $(N^{2n},h,J,\omega ,f_{N})$ be a complete gradient steady K\"{a}%
hler-Ricci soliton with exact K\"{a}hler form, and let $K$ be a compact
group of holomorphic isometries preserving $f_{N}$. Averaging a primitive
over $K$ gives a $K$-invariant one-form $\beta $ such that $d\beta =2\omega $%
. For $A\in K$ and $\ell >0$, define 
\begin{equation}
M_{A,\ell }=(N\times \mathbb{R})/\langle \sigma _{A,\ell }\rangle
\label{eq:steady-suspension}
\end{equation}%
with $\sigma _{A,\ell }(x,t)=(Ax,t+\ell ).$ The contactization tensors 
\begin{equation}
\eta =dt+\beta ,\qquad \xi =\partial _{t},\qquad g=\pi _{N}^{\ast }h+\eta
^{2},\qquad f=\pi _{N}^{\ast }f_{N},\qquad \Phi |_{\ker \eta }=-J^{\mathrm{%
hor}}  \label{eq:steady-suspension-tensors}
\end{equation}%
are invariant under $\sigma _{A,\ell }$. By Proposition~\ref%
{prop:regular-lift-data}, they define a complete Sasakian steady soliton on 
\begin{equation*}
N\times \mathbb{R}.
\end{equation*}%
The cyclic action is free because its nontrivial elements translate the $%
\mathbb{R}$-coordinate by a nonzero multiple of $\ell$, and it is properly
discontinuous for the same reason. Since it acts by isometries preserving
all of the tensors in \eqref{eq:steady-suspension-tensors}, the quotient $%
M_{A,\ell}$ is a complete Sasakian steady soliton.

\begin{theorem}
\label{thm:steady-suspension-classification} Let $M$ be a complete gradient
steady Sasaki-Ricci soliton with a complete smooth transverse development $%
(N^{2n},h,J,f_{N})$. Suppose 
\begin{equation*}
\func{Ric}_{h}>0
\end{equation*}
and $f_{N}$ has a critical point. Then $f_{N}$ has a unique critical point $%
o $, $N$ is biholomorphic to $\mathbb{C}^{n}$, and 
\begin{equation*}
K=\{B\in \func{Isom}(N,h):B_{\ast }J=JB_{\ast },\ f_{N}\circ B=f_{N}\}
\end{equation*}%
is compact and fixes $o$. After choosing a $K$-invariant primitive $d\beta
=2\omega $, the total soliton is isomorphic either to the line
contactization 
\begin{equation*}
N\times \mathbb{R}
\end{equation*}
or to a suspension $M_{A,\ell }$ with $A\in K$ and $\ell >0$. The line
contactization has a proper Reeb $\mathbb{R}$-action. For a suspension, the
Reeb types and periods are as follows.

\begin{enumerate}
\item If $A=\func{Id}$, the structure is regular with period $\ell $.

\item If $A$ has finite order $m>1$, the structure is quasi-regular but not
regular, its generic period is $m\ell $, and its effective orbifold quotient
is $[N/\langle A\rangle ]$. The leaf through $[(x,t)]$ has minimal period 
\begin{equation*}
j(x)\ell ,\qquad j(x)=\min \{j\geq 1:A^{j}x=x\}.
\end{equation*}%
In particular, the critical leaf over $o$ has period $\ell $.

\item If $A$ has infinite order, the structure is irregular. Put $C=%
\overline{\langle A\rangle }\subset K$. The closure of its Reeb flow is the
compact connected torus 
\begin{equation}
H=(C\times \mathbb{R})/\langle (A,\ell )\rangle ,\qquad \dim H=1+\dim
C^{0}\geq 2.  \label{eq:steady-suspension-torus}
\end{equation}%
A leaf is periodic exactly when $A^{j}x=x$ for some $j\geq 1$; the critical
leaf is always periodic and has period $\ell $.

\item For a fixed transverse soliton $(N,h,J,f_{N})$, two suspensions are
isomorphic preserving $(\xi ,\eta ,\Phi ,g,f)$ precisely when their positive
translation lengths agree and their monodromies are conjugate in $K$.
\end{enumerate}
\end{theorem}

\begin{proof}
The convexity argument in Proposition~\ref{prop:potential} gives uniqueness
of the critical point, properness of the increasing potential and spherical
positive levels. Morse coordinates near the minimum and the level flow show
that $N\cong\mathbb{R}^{2n}$. Conservation shows that $R_h$ attains its
maximum at $o$, so \cite[Theorem~1]{ChauTam2005} further gives that $N$ is
biholomorphic to $\mathbb{C}^n$. Every element of $K$ fixes $o$, so $K$ is a
closed subgroup of the compact isotropy group at $o$. Since $N$ is
contractible, $\omega$ is exact; averaging a primitive over $K$ gives $%
d\beta=2\omega$ with $\beta$ invariant.

Pass to the Riemannian universal cover $\widetilde M$ of $M$. The induced
development is locally trivial. Indeed, horizontal lifts of finite-length
base paths have the same length, and completeness prevents finite-time
escape. Transport along short radial geodesics therefore gives local
trivializations. By the homotopy sequence and contractibility of $N$, the
fibers in $\widetilde M$ are connected and simply connected. Thus they are
Reeb lines, and the development is a principal $\mathbb{R}$-bundle.
Trivializing it gives 
\begin{equation*}
\widetilde M=N\times\mathbb{R},\qquad \widetilde\eta=dt+\beta_0, \qquad
d\beta_0=2\omega.
\end{equation*}
The closed form $\beta_0-\beta$ is exact, so a change of the fiber
coordinate puts the contact form in the form $dt+\beta$.

We now determine the deck group $\Gamma$. As the deck group of a universal
cover, $\Gamma$ acts freely and properly discontinuously. A deck
transformation preserves the lifted Reeb foliation and induces a
potential-preserving holomorphic isometry $B\in K$ on $N$. Since it commutes
with $\partial_t$, it has the form $(x,t)\mapsto(Bx,t+c(x))$. Preservation
of $dt+\beta$ gives 
\begin{equation*}
dc=\beta-B^*\beta=0.
\end{equation*}
Thus $c$ is constant and $\Gamma$ is a subgroup of $K\times\mathbb{R}$. The
translation projection $p_{\mathbb{R}}:\Gamma\to\mathbb{R}$ is injective: if 
$p_{\mathbb{R}}(B,c)=0$, then $(B,0)$ fixes $(o,t)$ because every $B\in K$
fixes the unique critical point $o$; freeness forces $(B,0)$ to be the
identity.

The point $(o,0)$ has trivial stabilizer in $\Gamma$ because the deck action
is free. Proper discontinuity therefore implies that its orbit $%
\Gamma\cdot(o,0)$ is discrete. The image $p_{\mathbb{R}}(\Gamma)$ is
discrete. Otherwise there would be distinct $\gamma_j=(B_j,c_j)$ with $%
c_j\to0$. Since every $B_j$ fixes $o$, we have
$\gamma_j(o,0)=(o,c_j)\to(o,0)$, contradicting proper discontinuity.
Every discrete subgroup of $(\mathbb{R}%
,+)$ is either $\{0\}$ or $\ell\mathbb{Z}$ for a unique $\ell>0$.
Injectivity of $p_{\mathbb{R}}$ therefore implies that $\Gamma$ is trivial
or infinite cyclic, generated by a unique element $(A,\ell)$ with positive
translation length.

For a positive time $s$, the closing condition is 
\begin{equation*}
\varphi_s[(x,t)]=[(x,t)] \quad\Longleftrightarrow\quad s=j\ell,\quad A^jx=x
\quad\text{for some }j\in\mathbb{Z}_{>0}.
\end{equation*}
If $A$ has order $m$, its finite cyclic action is free on an open dense set,
giving generic period $m\ell$. The fixed point $o$ gives period $\ell$ for
the critical leaf. If $A$ has infinite order, the Reeb closure is $(C\times%
\mathbb{R})/\langle(A,\ell)\rangle$. This group is compact and abelian. It
is connected because $C/C^0$ is generated by the component of $A$, and is
therefore a torus of dimension $1+\dim C^0$.

Finally, an isomorphism preserving the soliton tensors lifts to $(B,c)$ on $%
N\times \mathbb{R}$. Conjugation sends the positive deck generator $(A,\ell
) $ to $(BAB^{-1},\ell )$. Conversely, such a conjugation descends to an
isomorphism. This proves the last assertion.
\end{proof}

\begin{remark}
Under $B^{T}\geq 0$ and $\func{Ric}^{T}>0$, Theorem~\ref%
{thm:positive-ricci-critical} supplies a critical point. If the Reeb closure
is compact, Proposition~\ref{prop:critical-leaf-development} constructs the
required development. If instead $H\cong\mathbb{R}$, the Reeb action is free
and proper, and its quotient $Q=M/\mathbb{R}$ is a complete smooth K\"ahler
manifold. The descended potential $f_Q$ has a critical point and $\func{Ric}%
_Q>0$. The convexity and level-flow argument of Proposition~\ref%
{prop:potential}, which needs no decay once a critical point exists, gives $%
Q\cong\mathbb{R}^{2n}$. Thus $Q$ is already simply connected and is itself
the complete transverse development. Proposition~\ref{prop:regular-lift-data}
gives $M$ as a line contactization $Q\times\mathbb{R}$. The suspension
theorem therefore proves Theorem~\ref{thm:intro-positive} in every
dimension. In dimension five, the developed model is either a product of two
scaled cigars or an irreducible K\"ahler steady surface with positive Ricci
curvature.
\end{remark}

\begin{proof}[Proof of Corollary~\protect\ref{cor:intro-decay-total}]
If $N$ is flat, choose the primitive
$\beta_H=\operatorname{Im}\langle z,dz\rangle$ of twice the Euclidean
K\"{a}hler form. The primitive in Theorem~\ref{thm:intro-decay} differs from
$\beta_H$ by an exact form, which can be absorbed in the fiber coordinate.
This gives the Heisenberg contactization with affine potential $f_H$.

Each deck transformation preserves the lifted tensors and commutes with the
Reeb flow. It therefore covers a holomorphic Euclidean isometry
$B:z\mapsto Az+v$ and has the form $(z,t)\mapsto(Az+v,t+\chi(z))$.
Preserving $\eta_H$ gives
\begin{equation*}
d\chi=\beta_H-B^*\beta_H
      =-d\operatorname{Im}\langle v,Az\rangle,
\end{equation*}
Thus the deck maps have the form \eqref{eq:flat-decay-deck}.
Preservation of $f_H$ is equivalent to $Aa=a$ and
$\operatorname{Re}\langle a,v\rangle=0$. Conversely, such transformations
preserve the Sasakian tensors and $f_H$, and a free, properly discontinuous
action gives a complete quotient with $R^T=0$.

If $N$ is nonflat, then $\func{Ric}_h>0$ and $f_N$ has a critical point by
Theorem~\ref{thm:intro-decay}. Theorem~\ref{thm:steady-suspension-classification}
gives the line and suspension models, the Reeb periods and flow closure, and
the isomorphism criterion. The fundamental groups follow because
$N\times\mathbb{R}$ is the universal cover in both models. Conversely, the
line contactization is complete by Proposition~\ref{prop:regular-lift-data},
and the suspension is its quotient by a free, properly discontinuous
cyclic group of soliton isometries. The transverse sectional curvature and
scalar curvature are those of $h$. The prescribed upper decay, or the
exponential decay for a scaled cigar, gives \eqref{eq:lifted-leaf-decay} by
\eqref{eq:lifted-leaf-distance}.
\end{proof}

In dimension three the nonflat base is a scaled cigar, and its
potential-preserving holomorphic isometries are the rotations about the tip.
Thus a nonflat suspension is generated by
$(z,t)\mapsto(e^{i\theta}z,t+\ell)$. Its Reeb flow is regular for the identity
rotation, quasi-regular for a nonidentity rotation of finite order, and
irregular when $\theta/(2\pi)$ is irrational.

\begin{definition}
\label{def:steady-soliton-approximation} A pointed complete gradient steady
Sasaki-Ricci soliton admits a \emph{quasi-regular steady soliton
approximation} if it is a pointed $C^{\infty }$ limit of complete
quasi-regular gradient steady Sasaki-Ricci solitons, including convergence
of the Sasakian tensors and of the potentials after additive normalization.
Such an approximation is \emph{developable} if all the approximating
orbifold quotients are developable.
\end{definition}

For suspensions, the approximation is obtained by varying the monodromy
within a compact group while keeping the transverse metric and potential
fixed.

\begin{proposition}
\label{prop:steady-monodromy-approximation} For any suspension %
\eqref{eq:steady-suspension} with $A$ of infinite order, there are
finite-order $A_{i}\in \overline{\langle A\rangle }$ and diffeomorphisms 
\begin{equation*}
\Psi _{i}:M_{A,\ell }\rightarrow M_{A_{i},\ell }
\end{equation*}
such that the pulled-back Sasakian tensors and potentials converge smoothly
on compact subsets to those of $M_{A,\ell }$.

\begin{enumerate}
\item If $m_{i}$ is the order of $A_{i}$, their generic Reeb periods satisfy 
$m_{i}\ell \rightarrow \infty $. The approximating orbifold quotients $%
[N/\langle A_{i}\rangle ]$ are developable. If $A$ and all $A_{i}$ fix a
point $o$, the leaf over $o$ retains period $\ell $ throughout the sequence.

\item For the same fixed transverse soliton, the regular contactizations 
\begin{equation*}
M_{\func{Id},L}=N\times (\mathbb{R}/L\mathbb{Z})
\end{equation*}
converge in the pointed smooth sense to the line contactization 
\begin{equation*}
N\times \mathbb{R}
\end{equation*}
as $L\rightarrow \infty $. Both assertions hold for nonflat $N$.
\end{enumerate}
\end{proposition}

\begin{proof}
Torsion elements are dense in every component of the compact abelian group $%
C=\overline{\langle A\rangle }$. Choose finite-order $A_{i}\rightarrow A$ in
the component of $A$, and write 
\begin{equation*}
A_{i}=A\exp B_{i},\qquad B_{i}\in \func{Lie}(C^{0}),\qquad B_{i}\rightarrow
0.
\end{equation*}%
Commutativity gives 
\begin{equation*}
\widetilde{\Psi }_{i}(x,t)=(\exp (tB_{i}/\ell )x,t),\qquad \widetilde{\Psi }%
_{i}\sigma _{A,\ell }=\sigma _{A_{i},\ell }\widetilde{\Psi }_{i}.
\end{equation*}%
Thus $\widetilde{\Psi }_{i}$ descends to a diffeomorphism $\Psi
_{i}:M_{A,\ell }\rightarrow M_{A_{i},\ell }$. Let $Y_{i}$ be the Killing
field generated by $B_{i}$. Pulling the target tensors back to the fixed
manifold $M_{A,\ell }$ gives 
\begin{equation}
\begin{aligned} \eta_i&=\beta+(1+\ell^{-1}\beta(Y_i))dt,\qquad
\xi_i=\partial_t-\ell^{-1}Y_i,\qquad f_i=f_N,\\ g_i((v,a),(w,b))
&=h(v+\ell^{-1}aY_i,w+\ell^{-1}bY_i) +\eta_i(v,a)\eta_i(w,b). \end{aligned}
\label{eq:steady-suspension-approximation}
\end{equation}%
Since $B_{i}\rightarrow 0$, these formulas give smooth convergence on
compact subsets, as well as convergence of $\Phi _{i}=-\nabla ^{g_{i}}\xi
_{i}$. Each structure is the pullback of a complete steady suspension. If
the orders $m_{i}$ had a bounded subsequence, a further subsequence would
satisfy $A_{i}^{m}=\func{Id}$ for one fixed $m$, forcing $A^{m}=\func{Id}$.
Hence $m_{i}\rightarrow \infty $. The period assertions follow from the
closing condition.

For the line-contactization limit, fix a relatively compact subset of $%
N\times \mathbb{R}$. Once $L$ exceeds the width of its $t$-coordinates, the
projection to $M_{\func{Id},L}$ is an embedding on this set and preserves
every tensor. Exhausting $N\times \mathbb{R}$ by such sets proves the
pointed smooth convergence.
\end{proof}

\begin{remark}
\label{rem:approx-transfer} The pulled-back structures in %
\eqref{eq:steady-suspension-approximation} generally have different CR
distributions. Type-I rational approximation preserves the CR structure, but
need not preserve the steady equation. The suspension approximation above
instead varies the monodromy while keeping the transverse steady soliton
fixed.
\end{remark}

\subsection{Decay from a Compact Reeb Orbit}

\label{subsec:compact-orbit-decay}

\begin{corollary}
\label{cor:irregular-decay-classification} Let $M^{5}$ be a complete
gradient steady Sasaki-Ricci soliton with $B^{T}\geq 0$ and $\func{Ric}%
^{T}>0 $. Suppose that the Reeb closure $H$ is compact and that, for one of
its orbits $Hp_{0}$, 
\begin{equation}
R^{T}(x)\longrightarrow 0\quad \text{uniformly as }d_{g}(Hp_{0},x)%
\longrightarrow \infty .  \label{eq:ambient-orbit-decay}
\end{equation}%
Then

\begin{enumerate}
\item $M$ is a suspension $M_{A,\ell }$ as in Theorem~\ref%
{thm:steady-suspension-classification}, over a complete simply connected
irreducible K\"{a}hler steady surface $N$ with $\func{Ric}_{h}>0$ and $%
B_{h}\geq 0$. The model $N$ is diffeomorphic to $\mathbb{R}^{4}$ and its
scalar curvature tends uniformly to zero at infinity.

\item The monodromy $A$ has finite order in the quasi-regular case and
infinite order in the irregular case. In the latter case $M$ has a
developable quasi-regular steady soliton approximation as in Definition~\ref%
{def:steady-soliton-approximation}.
\end{enumerate}
\end{corollary}

\begin{proof}
Theorem~\ref{thm:positive-ricci-critical} and Proposition~\ref%
{prop:critical-leaf-development} give a complete transverse development. The
suspension classification applies, and compactness of $H$ excludes the
line-contactization alternative.

Let $o$ be the unique critical point on $N$ and let $L$ be the critical Reeb
circle. We compare distance from $L$ with distance from $o$. Lift a curve
from $L$ to $[(x,t)]$ to $N\times \mathbb{R}$ and project it to $N$.
Projection decreases length, and the possible endpoints differ from $x$ by
powers of $A$, all of which fix $o$. This gives the lower bound in 
\begin{equation}
d_{g}(L,[(x,t)])=d_{h}(o,x),\qquad R^{T}([(x,t)])=R_{h}(x).
\label{eq:critical-circle-distance}
\end{equation}%
For the reverse bound, horizontally lift a minimizing geodesic from $o$ to $%
x $, choosing its initial fiber coordinate to end at $(x,t)$. The scalar
identity follows from the local transverse metric.

Since $L$ and $Hp_0$ are compact, 
\begin{equation*}
|d_g(L,p)-d_g(Hp_0,p)|\leq d_{\mathrm{H}}(L,Hp_0)<\infty,
\end{equation*}
where $d_{\mathrm{H}}$ denotes Hausdorff distance. Hence %
\eqref{eq:ambient-orbit-decay} transfers to uniform scalar decay on $N$.
Apply Theorem~\ref{thm:kahler-splitting}. Positive Ricci curvature excludes
flat factors, compact steady factors are ruled out by the trace equation,
and Proposition~\ref{prop:product-obstruction} excludes more than one
noncompact factor. Thus $N$ is irreducible. The assertions about monodromy
and approximation follow from Theorem~\ref%
{thm:steady-suspension-classification} and Proposition~\ref%
{prop:steady-monodromy-approximation}.
\end{proof}

\section{Explicit Examples and Constructions}

\label{sec:examples}

We first write the contactization data for the flat, product and irreducible
models. We then vary the monodromy and bundle data to illustrate the global
alternatives. The final two subsections describe integral circle lifts,
orbifold lifts, and irregular suspensions. Volume growth is treated in
Section~\ref{sec:asymptotic}.

\subsection{Contactizations and Cigar Products}

\label{subsec:contactizations}

We first introduce the method of contactization of an exact Kähler
soliton.

\begin{proposition}
\label{prop:contactization} Let $(N^{2n},h,J,\omega ,f_{N})$ be a complete
gradient steady K\"{a}hler-Ricci soliton, where $\omega (X,Y)=h(JX,Y)$ is
exact. Choose a one-form $A$ with $\mathrm{d}A=2\omega $. On $M=N\times 
\mathbb{R}$, set 
\begin{equation}
\eta =\mathrm{d}t+A,\qquad \xi =\partial _{t},\qquad g=\pi ^{\ast }h+\eta
^{2},\qquad f=\pi ^{\ast }f_{N}.  \label{eq:contactization}
\end{equation}%
On the horizontal distribution define $\Phi =-J^{\mathrm{hor}}$ and set $%
\Phi \xi =0$. Then $M$ is a complete gradient steady Sasaki-Ricci soliton.
The same construction on $N\times (\mathbb{R}/L\mathbb{Z})$ with $L>0$, is
complete and has a free Reeb circle action. Its transverse curvature is
precisely the curvature of $h$.
\end{proposition}

\begin{proof}
These are the trivial line and circle lifts in Proposition~\ref%
{prop:regular-lift-data}. Indeed, $d\eta =2\pi ^{\ast }\omega $, and its
Koszul calculation gives the Sasakian identity with $\Phi =-J^{\mathrm{hor}}$%
. The transverse connection is the pullback of the connection of $h$, so 
\begin{equation*}
\func{Ric}^{T}+\func{Hess}^{T}f=\pi ^{\ast }(\func{Ric}_{h}+\func{Hess}%
_{h}f_{N})=0.
\end{equation*}%
Completeness follows from the same proposition. The circle construction is
also the quotient of the line construction by the isometric translations $%
t\mapsto t+L$.
\end{proof}

\begin{proposition}
\label{prop:explicit-primitive} Let $\omega $ be a K\"{a}hler form on $%
\mathbb{C}^{n}$. In real linear coordinates, a primitive of $2\omega $ is
given by 
\begin{equation}
\beta _{x}(v)=2\int_{0}^{1}s\,\omega _{sx}(x,v)\,ds.
\label{eq:radial-primitive}
\end{equation}%
If a compact linear group preserves $\omega $, then this $\beta $ is
invariant under that group. If $\omega $ has a global K\"{a}hler potential
normalized by $\omega =(i/2)\partial \bar{\partial}u$, one may instead take 
\begin{equation}
\beta =\tfrac{1}{2}d^{c}u,\qquad d^{c}u=-du\circ J,\qquad dd^{c}u=4\omega .
\label{eq:potential-primitive}
\end{equation}
\end{proposition}

\begin{proof}
Writing $\omega =\tfrac{1}{2}\omega _{ij}dx^{i}\wedge dx^{j}$, the identity $%
d\omega =0$ and differentiation of \eqref{eq:radial-primitive} give 
\begin{equation*}
(d\beta )_{ij}=2\int_{0}^{1}\frac{d}{ds}\bigl(s^{2}\omega _{ij}(sx)\bigr)%
\,ds=2\omega _{ij}(x).
\end{equation*}%
A linear symmetry preserving $\omega $ commutes with radial contraction, so
it preserves $\beta $. For the second formula, 
\begin{equation*}
d^{c}=i(\bar{\partial}-\partial ),\qquad dd^{c}u=2i\partial \bar{\partial}%
u=4\omega .
\end{equation*}
\end{proof}

\subsection{Irreducible Steady K\"{a}hler-Ricci Soliton Models}

\label{subsec:irreducible-examples}

We construct complete regular five-dimensional gradient steady Sasaki--Ricci
solitons whose transverse model is not Cao's $U(2)$-invariant K\"{a}hler--Ricci
soliton.

\begin{proposition}
\label{prop:ac} There are complete regular five-dimensional gradient steady
Sasaki-Ricci solitons with positive transverse bisectional curvature whose
complete transverse model is irreducible and is not Cao's $U(2)$-invariant
soliton. Some such examples satisfy two-sided linear scalar curvature decay
on that model.
\end{proposition}

\begin{proof}
Choose a metric $h_{a}$ on $\mathbb{C}^{2}$, with parameter $a>0$, from
\cite[Corollary~1.3]{ApostolovCifarelli}, with the normalization of
\cite[equation~(6.1)]{ApostolovCifarelli}. It has positive sectional curvature
and is not Cao's metric. The K\"{a}hler identity 
\begin{equation*}
\func{Rm}(X,JX,JY,Y)=\func{Rm}(X,Y,Y,X)+\func{Rm}(X,JY,JY,X)
\end{equation*}%
gives positive bisectional curvature. Positive sectional curvature also
excludes a product decomposition. Since $H_{\mathrm{dR}}^{2}(\mathbb{C}%
^{2})=0$, Proposition~\ref{prop:contactization} gives regular Sasakian lifts.

In the variables $u=\xi _{1}$, $v=\xi _{2}$ of that reference, 
\begin{equation}
R_{h_{a}}=4a^{2}\frac{e^{-2au}-e^{-2av}}{v-u},\qquad 0\leq u\leq 1\leq
v,\qquad \rho _{h_{a}}\asymp v\quad (v\rightarrow \infty ).
\label{eq:ac-scalar}
\end{equation}%
At $u=v=1$, the expression is understood by continuity, with value $%
8a^3e^{-2a}$. See \cite[Appendix~A.2 and Lemma~5.8, equation~(5.18)]%
{ApostolovCifarelli}. For $v\geq 2$ and $e^{-2av}\leq e^{-2a}/2$, one has $%
e^{-2au}\leq1$ and $v-u\geq v-1\geq v/2$, hence 
\begin{equation*}
R_{h_a}\leq 4a^2\frac{1}{v-u}\leq\frac{8a^2}{v}.
\end{equation*}
On the other hand $u\leq1$ gives $e^{-2au}\geq e^{-2a}$, while the assumed
bound on $e^{-2av}$ gives $e^{-2au}-e^{-2av}\geq e^{-2a}/2$; since $v-u\leq
v $, this yields 
\begin{equation*}
\frac{2a^{2}e^{-2a}}{v}\leq R_{h_{a}}\leq \frac{8a^{2}}{v}.
\end{equation*}%
Here $\rho _{h_{a}}$ is distance from a fixed reference point. Its
comparison with $v$ gives the asserted two-sided linear bounds.
\end{proof}

\subsection{Elementary Steady Models}

Put $z_{j}=x_{j}+iy_{j}$, $h_{0}=\sum_{j=1}^{2}(\mathrm{d}x_{j}^{2}+\mathrm{d%
}y_{j}^{2})$, $\omega _{0}=\sum_{j}\mathrm{d}x_{j}\wedge \mathrm{d}y_{j}$,
and $f_{0}=0$.

For any $L>0,$ let $t\in \mathbb{R}/L\mathbb{Z}$ and set 
\begin{equation}
\eta _{0}=\mathrm{d}t+\sum_{j=1}^{2}(x_{j}\mathrm{d}y_{j}-y_{j}\mathrm{d}%
x_{j}),\quad g_{0}=h_{0}+\eta _{0}^{2}.  \label{eq:flat}
\end{equation}%
Then $\mathrm{d}\eta _{0}=2\omega _{0}$, and Proposition \ref%
{prop:contactization} applies. Since $H^{2}(\mathbb{C}^{2};\mathbb{Z})=0$,
the bundle is trivial: 
\begin{equation*}
M_{0}\cong \mathbb{C}^{2}\times S^{1},\quad e=0,\quad \pi _{1}(M_{0})=%
\mathbb{Z},\quad H^{2}(M_{0};\mathbb{Z})=0.
\end{equation*}%
The transverse geometry is flat, but the connection metric is not flat.

\begin{example}
\textbf{The cigar times a flat line : }
\end{example}

Write $z=x+iy=re^{i\theta }$ and $w=u+iv$. Consider 
\begin{equation}
\begin{split}
h_{c\times \mathbb{C}}& =\frac{\mathrm{d}x^{2}+\mathrm{d}y^{2}}{1+r^{2}}+%
\mathrm{d}u^{2}+\mathrm{d}v^{2}, \\
\omega _{c\times \mathbb{C}}& =\frac{\mathrm{d}x\wedge \mathrm{d}y}{1+r^{2}}+%
\mathrm{d}u\wedge \mathrm{d}v, \\
f_{c\times \mathbb{C}}& =-\log (1+r^{2}).
\end{split}
\label{eq:cigar}
\end{equation}%
Hamilton's cigar satisfies $\func{Ric}+\func{Hess}f=0$ in the displayed
scale. Its Gaussian curvature is $K=2/(1+r^{2})$ and its scalar curvature is 
$4/(1+r^{2})$. Define the globally smooth one-form 
\begin{equation}
A_{c}=\log (1+r^{2})\,\mathrm{d}\theta +u\mathrm{d}v-v\mathrm{d}u.
\label{eq:cigarA}
\end{equation}%
Indeed $\log (1+r^{2})\mathrm{d}\theta =(\log (1+r^{2})/r^{2})(x\mathrm{d}y-y%
\mathrm{d}x)$ extends smoothly at $r=0$, and $\mathrm{d}A_{c}=2\omega
_{c\times \mathbb{C}}$. Thus 
\begin{equation*}
M_{c\times \mathbb{C}}=\mathbb{C}^{2}\times (\mathbb{R}/L\mathbb{Z}),\quad
\eta =\mathrm{d}t+A_{c},\quad g=\pi ^{\ast }h_{c\times \mathbb{C}}+\eta ^{2}
\end{equation*}%
is complete and steady for every $L>0$. It has $e=0$, $\pi _{1}=\mathbb{Z}$,
and $R^{T}=4/(1+r^{2})$. Its potential is nonconstant.

\begin{example}
\textbf{The double cigar :}
\end{example}

On $Z=\mathbb{C}^{2}$ set 
\begin{align*}
h_{cc}& =\sum_{j=1}^{2}\frac{\mathrm{d}r_{j}^{2}+r_{j}^{2}\mathrm{d}\theta
_{j}^{2}}{1+r_{j}^{2}}, & f_{cc}& =-\sum_{j=1}^{2}\log (1+r_{j}^{2}), \\
\eta _{cc}& =\mathrm{d}t+\sum_{j=1}^{2}\log (1+r_{j}^{2})\mathrm{d}\theta
_{j}. & &
\end{align*}%
Again $\mathrm{d}\eta _{cc}=2\omega _{cc}$ and the angular one-forms extend
smoothly across the origins. The total space is $\mathbb{C}^{2}\times S^{1}$
with arbitrary Reeb period and 
\begin{equation*}
R^{T}=\frac{4}{1+r_{1}^{2}}+\frac{4}{1+r_{2}^{2}}.
\end{equation*}%
This is nonflat with two independent steady factors.

\subsection{The Integral Regular Lifting Theorem}

Let $(Z^4,h,J,\omega)$ be a K\"ahler surface with $\omega(X,Y)=h(JX,Y)$. We
take the circle action of a Reeb field $\xi$ to have minimal period $L>0$;
thus $\eta(\xi)=1$ and $\mathrm{d}\eta=2\,\pi^*\omega$, where $\pi^*$
denotes pullback along the bundle projection $\pi\colon M\to Z$ (not
multiplication by $2\pi$).

\begin{theorem}
\label{thm:lift} Suppose $(Z,h,J,f_{Z})$ is a connected complete gradient
steady K\"{a}hler--Ricci soliton: 
\begin{equation*}
\func{Ric}_{h}+\func{Hess}_{h}f_{Z}=0,\qquad \func{Hess}_{h}^{2,0}f_{Z}=0.
\end{equation*}%
Fix $L>0$. A principal Reeb-circle connection lift of period $L$ with $%
\mathrm{d}\eta =2\,\pi ^{\ast }\omega $ exists if and only if 
\begin{equation}
\frac{2}{L}[\omega ]\in \func{im}\bigl(H^{2}(Z;\mathbb{Z})\longrightarrow
H^{2}(Z;\mathbb{R})\bigr).  \label{eq:int}
\end{equation}%
For each chosen integral lift $e$ of this real class there is a principal
circle bundle 
\begin{equation*}
\pi :M\rightarrow Z
\end{equation*}
with Euler class $e$, a connection form $\eta $, and a complete Sasakian
metric 
\begin{equation}
g=\pi ^{\ast }h+\eta ^{2},\qquad f=\pi ^{\ast }f_{Z},  \label{eq:metric}
\end{equation}%
which solves \eqref{eq:steady}. Conversely every regular Sasakian steady
gradient soliton with this quotient and period satisfies \eqref{eq:int}. For
fixed bundle, curvature and period, connections differ by closed one-forms
on $Z$. If $H^{1}(Z;\mathbb{R})=0,$ they are gauge-equivalent by a bundle
automorphism covering the identity.
\end{theorem}

\begin{proof}
Let $X=(L/2\pi)\xi$ be the fundamental field of the standard $2\pi$-periodic
circle action, with connection form $\alpha=(2\pi/L)\eta$. Writing $%
d\alpha=\pi^*\Omega_\alpha$, Chern--Weil theory gives 
\begin{equation*}
e_{\mathbb{R}}=[\Omega_\alpha/(2\pi)]=\frac{2}{L}[\omega].
\end{equation*}
Conversely, for each chosen integral lift $e$, the construction in
Proposition~\ref{prop:regular-lift-data} gives the required principal circle
bundle and connection. The same proposition proves the Sasakian identity,
the steady equation and completeness of the connection metric.

Finally, two connections on the fixed bundle with the same curvature satisfy 
$\eta_2-\eta_1=\pi^*a$ for a closed one-form $a$ on $Z$. If $H^1(Z;\mathbb{R}%
)=0$, write $a=du$. The bundle automorphism $\Psi(p)=\varphi_{-u(\pi(p))}(p)$%
, where $\varphi_t$ denotes fiber translation with period $L$, satisfies $%
\Psi^*\eta_2=\eta_1$. This proves the claimed gauge equivalence.
\end{proof}

\begin{example}
\textbf{Cigar Times an Elliptic Curve :}
\end{example}

Let $E=\mathbb{C}/\Lambda $ be an elliptic curve with flat K\"{a}hler metric 
$h_{E}$, form $\omega _{E}$, and area $A_{E}=\int_{E}\omega _{E}>0$. Put 
\begin{equation*}
Z=\mathbb{C}\times E,\quad h=h_{c}\oplus h_{E},\quad \omega =\omega
_{c}+\omega _{E},\quad f_{Z}=-\log (1+|z|^{2}).
\end{equation*}%
Since $Z$ retracts onto $E$, $H^{2}(Z;\mathbb{Z})=\mathbb{Z}$. Choose $d\geq
1$ and let 
\begin{equation*}
P_{d}\rightarrow E
\end{equation*}%
have Euler number $d$. By Theorem \ref{thm:lift}, the necessary and
sufficient period is 
\begin{equation}
L=\frac{2A_{E}}{d}.  \label{eq:ellperiod}
\end{equation}%
Let $\eta _{E}$ be the real-valued connection on $P_{d}$ with $\eta _{E}(\xi
)=1$, Reeb period $L$, and $\mathrm{d}\eta _{E}=2\func{pr}_{E}^{\ast }\omega
_{E}$. On 
\begin{equation*}
M=\mathbb{C}\times P_{d},
\end{equation*}
set 
\begin{equation}
\eta =\eta _{E}+\log (1+r^{2})\mathrm{d}\theta ,\qquad g=\pi ^{\ast
}(h_{c}+h_{E})+\eta ^{2}.  \label{eq:ellmetric}
\end{equation}%
The curvature is $2(\omega _{c}+\omega _{E})$, so \eqref{eq:ellmetric} is a
complete nonconstant-potential steady soliton.

\begin{proposition}
\label{prop:nil} The compact $P_d$ is an oriented Heisenberg nilmanifold.
With a choice of orientation and generators, 
\begin{equation*}
\pi_1(M)\cong\langle a,b,c\mid c\text{ central},\ [a,b]=c^d\rangle,\qquad
H^2(M;\mathbb{Z})\cong\mathbb{Z}^2\oplus\mathbb{Z}/d\mathbb{Z}.
\end{equation*}
\end{proposition}

\begin{proof}
The circle bundle over the torus has central fiber subgroup $\langle
c\rangle $ and extension class $d$, yielding the indicated presentation. The
Gysin sequence gives $H^2(P_d;\mathbb{Z})$ as an extension of $\mathbb{Z}^2$
by $\mathbb{Z}/d\mathbb{Z}$, which splits since $\mathbb{Z}^2$ is free.
Contractibility of $\mathbb{C}$ yields the assertions for $M$.
\end{proof}

\begin{example}
\textbf{Eguchi--Hanson transverse geometry :}
\end{example}

This example illustrates how nonflat Ricci-flat transverse geometry can
produce nontrivial circle topology.

\begin{proposition}
\label{prop:EH} Let $Q=\func{Tot}\mathcal{O}_{\mathbb{CP}^{1}}(-2)$ carry
the complete Eguchi--Hanson Ricci-flat K\"{a}hler metric \cite{EguchiHanson}%
, with zero section $E_{0}$ of area $A_{0}$. For $q\geq 1$, choose Euler
number $q$ on $E_{0}$ and period 
\begin{equation}  \label{eq:EHperiod}
L=\frac{2A_{0}}{q}.
\end{equation}
The regular lift $M_{q}$ is complete and steady with constant potential. Its
restriction over $E_{0}$ is $L(q,1)$, up to orientation. Since $Q$ retracts
onto $E_{0}$, 
\begin{equation*}
\pi _{1}(M_{q})\cong \mathbb{Z}/q\mathbb{Z},\qquad H^{2}(M_{q};\mathbb{Z}%
)\cong \mathbb{Z}/q\mathbb{Z}.
\end{equation*}%
In particular, $M_{1}$ is simply connected with nonflat Ricci-flat
transverse metric. It cannot satisfy $B^{T}\geq 0$, by Lemma~\ref%
{lem:curvature-control}.
\end{proposition}

\begin{proof}
The boundary map $\pi_2(Q)=\mathbb{Z}\to\pi_1(S^1)=\mathbb{Z}$ is
multiplication by $q$, and the long exact sequence gives $\pi_1(M_q)=\mathbb{%
Z}/q\mathbb{Z}$. The Gysin sequence gives $H^2(M_q;\mathbb{Z})=H^2(Q;\mathbb{%
Z})/\langle e\rangle=\mathbb{Z}/q\mathbb{Z}$.
\end{proof}

\begin{remark}
Eguchi--Hanson is steady with constant potential, not a nontrivial steady
Ricci soliton. Its value here is the nonflat Ricci-flat transverse geometry
and nontrivial circle topology.
\end{remark}

\subsection{Orbifold Lifts}

\label{subsec:bundle-examples}

\begin{proposition}
\label{prop:orbifold-integral-lift} Let $(\mathcal{Q},h,J,\omega ,f_{%
\mathcal{Q}})$ be a complete effective K\"{a}hler orbifold carrying a
gradient steady soliton. Suppose that an orbifold principal circle bundle $%
\pi :P\rightarrow \mathcal{Q}$ has integral Euler class $e_{\mathrm{orb}}$
satisfying 
\begin{equation}
(e_{\mathrm{orb}})_{\mathbb{R}}=\frac{2}{L}[\omega ]_{\mathrm{orb}},\qquad
L>0.  \label{eq:orbifold-integral-lift}
\end{equation}%
Here integral orbifold cohomology means classifying-space cohomology. Assume
that the fiber character is injective on every local point stabilizer. Then

\begin{enumerate}
\item $P$ is a complete gradient steady Sasaki-Ricci soliton and admits a
connection $\eta $ with $d\eta =2\pi ^{\ast }\omega $. Its metric is $%
g=\pi^{\ast}h+\eta^2$, its basic potential is $\pi^{\ast}f_{\mathcal{Q}}$,
and its generic Reeb period is $L$.

\item For local cyclic isotropy $\mathbb{Z}_{m}$ acting on the fiber with
weight $a$, the local total space is smooth precisely when $\gcd(a,m)=1$. In
that case the Reeb orbit over the point has minimal period $L/m$.
\end{enumerate}
\end{proposition}

\begin{proof}
Choose a $2\pi $-periodic connection $\alpha _{0}$ of curvature $\pi ^{\ast
}\Omega $. Chern--Weil and orbifold de Rham theory give 
\begin{equation*}
\left[ \frac{\Omega }{2\pi }\right] =\frac{2}{L}[\omega ],\qquad \Omega -%
\frac{4\pi }{L}\omega =d\beta
\end{equation*}%
for a global orbifold one-form $\beta $ \cite{Satake1956,BoyerGalicki}. Then 
\begin{equation*}
\eta =\frac{L}{2\pi }(\alpha _{0}-\pi ^{\ast }\beta ),\qquad d\eta =2\pi
^{\ast }\omega .
\end{equation*}%
The local calculation of Proposition~\ref{prop:regular-lift-data} proves the
Sasakian identity and the steady equation in every chart. A stabilizer
element fixes a point upstairs precisely when its fiber character is
trivial, proving the smoothness criterion. For $\mathbb{Z}_{m}$ this is
equivalent to $\gcd (a,m)=1$. In that case the image of the fiber character
is the full subgroup of order $m$ in $\mathbb{R}/L\mathbb{Z}$, so the fiber
over that point has minimal Reeb period $L/m$.

Finally, a finite-length curve projects to a convergent curve in the
complete orbifold base. A compact neighborhood of the limit has compact
inverse image, since the fibers are circles and the local groups are finite.
Hence the curve cannot escape to infinity, proving completeness.
\end{proof}

\begin{example}
Irregular nonflat globally developable suspensions
\end{example}

We now work with the contactization of \eqref{eq:cigar}: let $\widetilde{M}=%
\mathbb{C}^{2}\times \mathbb{R}$ with $\widetilde{\eta }=\mathrm{d}t+A_{c}$, 
$\widetilde{g}=h_{c\times \mathbb{C}}+\widetilde{\eta }^{2}$, $\widetilde{%
\xi }=\partial _{t}$, and $\widetilde{f}=-\log (1+|z|^{2})$.

\begin{theorem}
\label{thm:irr} Fix $\ell >0$ and angles $\alpha ,\beta \in \mathbb{R}$. The
diffeomorphism 
\begin{equation*}
\gamma (z,w,t)=(e^{i\alpha }z,e^{i\beta }w,t+\ell )
\end{equation*}%
generates a free properly discontinuous group of Sasakian soliton
automorphisms. Its quotient 
\begin{equation*}
M_{\alpha ,\beta ,\ell }=\widetilde{M}/\langle \gamma \rangle
\end{equation*}%
is a complete five-dimensional steady gradient Sasaki--Ricci soliton with
nonconstant potential and nonflat transverse curvature, with $\pi
_{1}(M_{\alpha ,\beta ,\ell })\cong \mathbb{Z}$. The lifted Reeb foliation
on its universal cover is the fiber foliation of the complete submersion 
\begin{equation*}
\widetilde{M}\rightarrow \mathbb{C}^{2}.
\end{equation*}
If at least one of $\alpha /2\pi $ or $\beta /2\pi $ is irrational, the
quotient Reeb foliation is irregular.
\end{theorem}

\begin{proof}
Both terms of $A_{c}$ are invariant under the corresponding planar
rotations. Hence $\gamma ^{\ast }\widetilde{\eta }=\widetilde{\eta }$, and $%
\gamma $ also preserves $h$, the transverse complex structure, and $%
\widetilde{f}$. Translation by $\ell $ in the $t$-coordinate makes the
cyclic action free and properly discontinuous. A quotient of a complete
Riemannian manifold by a freely acting discrete isometry group is complete.
The universal cover $\mathbb{C}^{2}\times \mathbb{R}$ is simply connected,
so the quotient fundamental group is $\mathbb{Z}$.

A Reeb orbit through $(z,w,t)$ closes on the quotient precisely when for
some nonzero integer $k$, 
\begin{equation*}
e^{ik\alpha }z=z,\qquad e^{ik\beta }w=w.
\end{equation*}%
If an irrational angle corresponds to a nonzero coordinate, no such $k$
exists. Thus the Reeb foliation is irregular whenever at least one angle is
irrational. The projection $(z,w,t)\mapsto (z,w)$ is a complete smooth
transverse development on the universal cover, as claimed.
\end{proof}

\begin{example}
\label{ex:nonflat-irregular-steady}

\begin{enumerate}
\item The flat, cigar-line, two-cigar, Cao and Apostolov--Cifarelli models
above have exact K\"{a}hler forms and potential-preserving rotational
symmetries. Infinite-order rotations give complete irregular suspensions on 
\begin{equation*}
\mathbb{R}^{4}\times \mathbb{S}^{1},
\end{equation*}%
and finite-order rotations give the quasi-regular steady approximations of
Proposition~\ref{prop:steady-monodromy-approximation}.

\item For the Cao and Apostolov--Cifarelli models, the common transverse
metric has positive bisectional curvature and two-sided linear scalar
curvature decay. Choosing a nonflat model and letting the regular period
tend to infinity instead yields its nonflat line contactization.
\end{enumerate}
\end{example}

\section{Asymptotic Geometry and Volume Growth}

\label{sec:asymptotic}

We conclude with volume-growth and asymptotic calculations for the explicit
models constructed above. These examples also clarify the distinction
between transverse volume growth and the volume growth of the total Sasakian
space, particularly for suspensions with compact Reeb fundamental domain.

\begin{lemma}
\label{lem:cigend} Set $s=\func{arsinh}r$. The cigar metric, potential, and
scalar curvature become 
\begin{equation*}
h_{c}=\mathrm{d}s^{2}+\tanh ^{2}s\,\mathrm{d}\theta ^{2},\quad f_{c}=-2\log
\cosh s,\quad R_{c}=4\func{sech}^{2}s.
\end{equation*}%
As $s\rightarrow \infty $ the metric tends exponentially in every derivative
to the flat cylinder $\mathrm{d}s^{2}+\mathrm{d}\theta ^{2}$ on unit-width
longitudinal windows. Its Gaussian curvature is $2\func{sech}^{2}s$.
Moreover $f_{c}=-2s+2\log 2+O(e^{-2s})$.
\end{lemma}

\begin{proof}
Since $r=\sinh s$, $\mathrm{d}r^{2}/(1+r^{2})=\mathrm{d}s^{2}$ and $%
r^{2}/(1+r^{2})=\tanh ^{2}s$. The remaining formulas follow by direct
substitution and Taylor expansion of $\tanh s$ and $\log \cosh s$.
\end{proof}

\begin{proposition}
\label{prop:growth} Let $V_{M}(R)$ be the volume of a geodesic ball of
radius $R$ in a regular circle connection lift, based at a fixed point.
There are positive constants $c,C$ (depending on the model and the base
point) such that, for all sufficiently large $R$, 
\begin{equation*}
cR^{d}\leq V_{M}(R)\leq CR^{d},
\end{equation*}%
with the following exponents: 
\begin{equation*}
\begin{array}{c|ccccc}
Z & \mathbb{C}^{2} & \text{cigar}\times \mathbb{C} & \text{cigar}\times 
\text{cigar} & \text{cigar}\times E & \text{Eguchi--Hanson} \\ \hline
d & 4 & 3 & 2 & 1 & 4%
\end{array}%
\end{equation*}%
The quasi-regular finite quotient of the cigar-line lift has exponent $3$.
The irrational cigar-line suspensions also have exponent $3$, as established
independently in Theorem~\ref{thm:irrvolume}. Suspensions over the
double-cigar model have exponent $2$ (Corollary~\ref%
{cor:double-cigar-suspension-growth}).
\end{proposition}

\begin{proof}
The cigar has linear area growth in geodesic distance: in the coordinate $s=%
\func{arcsinh}r$, the area of a geodesic disk of radius $R$ is 
\begin{equation*}
2\pi \int_{0}^{R}\tanh s\,\mathrm{d}s=2\pi \log \cosh R\sim 2\pi R.
\end{equation*}%
The Euclidean plane $\mathbb{C}$ has quadratic area growth $\pi R^{2}$. For
the product of two cigar surfaces the geodesic distance is $\sqrt{%
d_{1}^{2}+d_{2}^{2}}$, and the volume of a metric ball of radius $R$ is 
\begin{equation*}
\int_{0}^{R}\int_{0}^{\sqrt{R^{2}-s_{1}^{2}}}(2\pi \tanh s_{1})(2\pi \tanh
s_{2})\,\mathrm{d}s_{2}\,\mathrm{d}s_{1}.
\end{equation*}
For a large $R$ the integrand is approximately $(2\pi )^{2}$ on the
quarter-disk, giving growth $\sim (2\pi )^{2}\cdot \pi R^{2}/4=\pi ^{3}R^{2}$%
, i.e.\ exponent $2$. Products therefore have exponents $3$ (cigar times $%
\mathbb{C}$), $2$ (cigar times cigar), and $1$ when the second factor is
compact (cigar times $E$). Eguchi--Hanson is ALE of real dimension four and
has fourth-order volume growth. For a regular lift with fixed compact circle
fibers, projection is distance nonincreasing (since $|\mathrm{d}\pi
(v)|_{h}\leq |v|_{g}$ by the horizontal--vertical decomposition).
Conversely, horizontal lifting of a minimizing base path followed by a fiber
arc of length at most $L/2$ shows that the inverse image of a base ball of
radius $R-L/2$ lies in a total-space ball of radius $R$ (with
basepoint-dependent choice of initial fiber point). The volume form is the
base volume form times the fiber arclength, so the same polynomial growth
order follows. For the finite quotient, let $p:\widehat{M}\rightarrow M$ be
the finite Riemannian covering of degree $m$. Choosing lifts of a base
point, the distance from any lift to the chosen lift is bounded by a
constant depending only on this finite orbit. The lifted inverse image of a
ball is trapped between finite unions of balls of radii $R\pm C$, so finite
quotients preserve the polynomial growth exponent.
\end{proof}

\begin{lemma}
\label{lem:suspension-volume-comparison} Let $M_{A,\ell }$ be a connection
suspension as in \eqref{eq:steady-suspension}--%
\eqref{eq:steady-suspension-tensors}, with complete base $(N,h)$, and
suppose that $A$ fixes $o\in N$. Let $L=\{[o,t]\}$ be the Reeb circle over $%
o $. For $p_{0}\in L$ and $R>\ell /2$, one has 
\begin{equation}
\ell \,\func{Vol}_{h}B_{h}(o,R-\ell /2)\leq \func{Vol}_{g}B_{g}(p_{0},R)\leq
\ell \,\func{Vol}_{h}B_{h}(o,R).  \label{eq:suspension-volume-comparison}
\end{equation}%
In particular, any leading polynomial volume asymptotic of the base
transfers to the suspension, with coefficient multiplied by $\ell $.
\end{lemma}

\begin{proof}
Since $A$ fixes $o$, the function $r([x,t])=d_{h}(o,x)$ is well defined.
Projection of a lifted curve to $N$ does not increase length, and horizontal
lifting of a minimizing base geodesic, with prescribed terminal point, gives
the reverse inequality. Hence 
\begin{equation*}
d_{g}(L,[x,t])=r([x,t]).
\end{equation*}%
The connection volume form on $0\leq t<\ell $ is $dt\,d\func{vol}_{h}$, so
the radius-$R$ tube about $L$ has volume $\ell \,\func{Vol}_{h}B_{h}(o,R)$.
The circle $L$ has length $\ell $. For $p_{0}\in L$, a shortest arc in $L$
has length at most $\ell /2$, which proves %
\eqref{eq:suspension-volume-comparison} and the asymptotic conclusion.
\end{proof}

\begin{theorem}
\label{thm:irrvolume} For every $\ell >0$ and every $\alpha ,\beta \in 
\mathbb{R}$ (including irrational angles), the complete suspension $%
(M_{\alpha ,\beta ,\ell },g)$ of Theorem~\ref{thm:irr} satisfies 
\begin{equation*}
cR^{3}\leq \func{Vol}_{g}B_{g}(p,R)\leq CR^{3}\qquad (R\geq R_{0}),
\end{equation*}%
where $c,C,R_{0}>0$ may depend on $p$ and $\ell $ but the constants can be
chosen independently of $\alpha ,\beta $. In particular the volume-growth
exponent is exactly $3$, without any Diophantine hypothesis on the angles.
More precisely, 
\begin{equation*}
\func{Vol}_{g}B_{g}(p,R)\sim \frac{4\pi ^{2}\ell }{3}R^{3}.
\end{equation*}
\end{theorem}

\begin{proof}
For the present cigar-line model, take $o=(0,0)$ and $p_{0}=[0,0,0]$. In
cigar polar coordinates $z=\sinh s\,e^{i\theta }$ and Euclidean polar
coordinates $w=\rho e^{i\varphi }$, the base distance from $o$ is $\sqrt{%
s^{2}+\rho ^{2}}$ and the angularly integrated volume density is $(2\pi
\tanh s\,ds)(2\pi \rho \,d\rho )$. Consequently 
\begin{equation*}
\func{Vol}_{h}B_{h}(o,R)=2\pi ^{2}\int_{0}^{R}\tanh s\,(R^{2}-s^{2})\,ds\sim 
\frac{4\pi ^{2}}{3}R^{3}.
\end{equation*}%
Indeed, after setting $s=Rv$, dominated convergence gives $\int_{0}^{1}\tanh
(Rv)(1-v^{2})\,dv\longrightarrow 2/3$. Lemma~\ref%
{lem:suspension-volume-comparison} proves the asserted asymptotic and bounds
at $p_{0}$, independently of the rotation angles. For a fixed $p=[z,w,t]$,
put $D=d_{g}(p,p_{0})$. The same horizontal-lifting argument gives $D\leq
\ell /2+\sqrt{s^{2}+\rho ^{2}}$, a bound independent of the angles. Finally, 
$B(p_{0},R-D)\subset B(p,R)\subset B(p_{0},R+D)$ proves the statements at $p$%
.
\end{proof}

\begin{corollary}
\label{cor:double-cigar-suspension-growth} Let $M_{A,\ell }$ be any
suspension of Theorem~\ref{thm:steady-suspension-classification} whose
transverse development is a product of two scaled cigars and whose monodromy 
$A$ preserves the soliton potential. Then 
\begin{equation*}
cR^{2}\leq \func{Vol}B(p,R)\leq CR^{2}
\end{equation*}%
for all sufficiently large $R$. Thus the suspension has the same polynomial
volume-growth exponent $2$ as its double-cigar transverse model; the compact
fundamental interval in the Reeb coordinate does not add an extra power of $%
R $. If the transverse metric is $ag_{\mathrm{cig}}\times bg_{\mathrm{cig}}$%
, then more precisely 
\begin{equation*}
\func{Vol}B(p,R)\sim \ell \pi ^{3}\sqrt{ab}\,R^{2}.
\end{equation*}
\end{corollary}

\begin{proof}
Let $o$ be the pair of cigar tips. Since $o$ is the unique critical point of
the soliton potential, the potential-preserving monodromy fixes $o$. Use
geodesic radial coordinates $s_{1},s_{2}$ for $ag_{\mathrm{cig}}$ and $bg_{%
\mathrm{cig}}$. Their angularly integrated area densities are 
\begin{equation*}
2\pi \sqrt{a}\,\tanh (s_{1}/\sqrt{a})\,ds_{1}\text{ \ \ \textrm{and} \ \ \ }%
2\pi \sqrt{b}\,\tanh (s_{2}/\sqrt{b})\,ds_{2}.
\end{equation*}%
The base ball is the quarter-disk $s_{1}^{2}+s_{2}^{2}<R^{2}$ in these
coordinates. Rescaling by $R$ and applying dominated convergence therefore
gives 
\begin{equation*}
\func{Vol}_{h}B_{h}(o,R)\sim 4\pi ^{2}\sqrt{ab}\,\frac{\pi R^{2}}{4}=\pi ^{3}%
\sqrt{ab}\,R^{2}.
\end{equation*}%
Lemma~\ref{lem:suspension-volume-comparison} and a fixed change of basepoint
prove the conclusion. The argument is independent of whether $A$ has finite
or infinite order.
\end{proof}

\section*{Declaration of AI use}
The authors used ChatGPT (OpenAI) during the
preparation of the manuscript for language editing and exploratory calculations.
All mathematical arguments and computations were independently checked by the
authors.

\end{document}